\documentclass[a4paper,english,cleveref, autoref, thm-restate,
nolineno,
]{socg-lipics-v2021}

\pdfoutput=1 %uncomment to ensure pdflatex processing (mandatatory e.g. to submit to arXiv)
\hideLIPIcs  %uncomment to remove references to LIPIcs series (logo, DOI, ...), e.g. when preparing a pre-final version to be uploaded to arXiv or another public repository

\usepackage{bm}
\usepackage{tikz-cd}
\usepackage{mathtools}
\usepackage{numprint}
\npthousandsep{,}
\npthousandthpartsep{}
\usepackage{booktabs}
\usepackage{multirow}
\usepackage{cellspace}
\usepackage{stackengine}
\usepackage{diagbox}
\setstackEOL{\\}
\makeatletter
\newcommand\dashedmidrule[1]{%
  \noalign{\vskip\dimexpr\aboverulesep+\cmidrulewidth-\ht\@arstrutbox\relax}%
  \multicolumn{#1}{@{}c@{}}{\xleaders\hbox to 4pt{\hss\rule{2pt}{\cmidrulewidth}\hss}\hfill\kern0pt}\\
  \noalign{\vskip\dimexpr\belowrulesep-\dp\@arstrutbox\relax}}
\makeatother

\DeclareMathOperator{\Aut}{Aut}
\DeclareMathOperator{\End}{End}

\title{Fixed-parameter tractable computation of Reshetikhin--Turaev knot polynomials via Tensor Networks} 

\titlerunning{FPT computation of RT knot polynomials via Tensor Networks} %TODO optional, please use if title is longer than one line

\author{Shana Yunsheng Li}{Dept. of Math., University of Illinois Urbana-Champaign, USA \and \url{https://shana-y-li.github.io} }{yl202@illinois.edu}{https://orcid.org/0009-0006-5320-106X}{}%TODO mandatory, please use full name; only 1 author per \author macro; first two parameters are mandatory, other parameters can be empty. Please provide at least the name of the affiliation and the country. The full address is optional. Use additional curly braces to indicate the correct name splitting when the last name consists of multi$ple name parts.

\authorrunning{S. Y. Li} %TODO mandatory. First: Use abbreviated first/middle names. Second (only in severe cases): Use first author plus 'et al.'

\Copyright{Shana Yunsheng Li} %TODO mandatory, please use full first names. LIPIcs license is "CC-BY";  http://creativecommons.org/licenses/by/3.0/

%TODO mandatory: Please choose ACM 2012 classifications from https://dl.acm.org/ccs/ccs_flat.cfm 

 %TODO mandatory; please add comma-separated list of keywords

\category{} %optional, e.g. invited paper

\funding{This material is based upon work supported by the U.S. National Science Foundation under Grant No. DMS-2424139 while the author was in residence at the Simons Laufer Mathematical Sciences Institute in Berkeley,
California, during the Spring 2026 semester. The author was also supported by the U.S. National Science Foundation under Grant No. DMS-2303572 during the Spring and Summer 2026 semesters.}%optional, to capture a funding statement, which applies to all authors. Please enter author specific funding statements as fifth argument of the \author macro.

\acknowledgements{The author wishes to thank Nathan Dunfield, Stavros Garoufalidis, Saul Schleimer, and Jonathan Spreer for useful conversations. The author also specially thanks Mingde Ren, from whom the author learned about tensor networks. }%optional

\EventEditors{John Q. Open and Joan R. Access}
\EventNoEds{2}
\EventLongTitle{42nd Conference on Very Important Topics (CVIT 2016)}
\EventShortTitle{CVIT 2016}
\EventAcronym{CVIT}
\EventYear{2016}
\EventDate{December 24--27, 2016}
\EventLocation{Little Whinging, United Kingdom}
\EventLogo{}
\SeriesVolume{42}
\ArticleNo{23}
\begin{document}

\maketitle

\begin{abstract}
  We give a thorough analysis of the time complexity of computing Reshetikhin--Turaev knot polynomials via tensor contractions on the associated tensor networks, showing that the computation is fixed-parameter tractable with respect to a parameter at most linear in the tree-width of the input knot diagram. 
  In particular, this gives a simpler proof to the sub-exponential bound $e^{O(\sqrt{n})}$ for the time complexity of computing any Reshetikhin--Turaev knot polynomial. 
  We present along with this paper, in \texttt{SnapPy}, the first general-purpose implementation of this algorithm, which computes any Reshetikhin--Turaev knot polynomial given its $R$-matrix and ribbon element. We stress the remarkable efficiency of our implementation by demonstrating that: (1) our implementation is several orders of magnitude faster than existing packages on computing the colored Jones polynomial; (2) our implementation is able to compute the Jones polynomial of (i) census knots with large width and over 200 crossings and (ii) some generic random knots with over 2000 crossings, in several hours using less than 10GB RAM, whereas \texttt{Regina} is not able to compute the Jones polynomial of these knots within a memory constraint of 160GB. 
\end{abstract}

\section{Introduction}
\label{sec.intro}

The discovery of the Jones polynomial \cite{Jones} revealed deep connections of low-dimensional topology with physics \cite{Witten}, which gave birth to the field of quantum topology. Building upon this foundation, the Reshetikhin--Turaev functor \cite{Turaev:YB,RT:ribbon} was later established as a powerful and systematic approach to defining numerous knot polynomials, including the colored Jones polynomials \cite{KR:colored-jones,Turaev:YB}, the colored HOMFLY polynomials, and the colored Links--Gould polynomials \cite{GK:multi,Vn}.
As is the case for many invariants, efficient computation is essential for studying their properties \cite{Bar-Natan:Khovanov,Shumakovitch:odd,GL:patterns,GLY:LG}, and is thus naturally an important topic. However, it was later understood that  Reshetikhin--Turaev knot polynomials are generally $\#\mathrm{P}$-hard to compute, as it was proved by Jaeger, Vertigan and Welsh for the Jones polynomial \cite{JVW:complexity}, one of the simplest Reshetikhin--Turaev knot polynomials. $\#\mathrm{P}$-hard problems are analogues of $\mathrm{NP}$-hard problems, where one computes a value (in this case, the knot polynomial) rather than deciding a yes/no question. Thus for $\#\mathrm{P}$-hard problems, a natural goal is an algorithm that is \textit{fixed-parameter tractable} \cite{CFKLMPPS:parametrized}, which in general means that the input can be subdivided into families $\mathcal{F}_i$ parametrized by $i\in \mathbb{N}$, such that the time complexity for an input in $\mathcal{F}_i$ of size $n$ is 
\begin{equation}
  \label{eq.FPT}
  O(f(i)\cdot p(n)),
\end{equation}
where $p$ is a polynomial and $f$ is an arbitrary computable function. We say that an algorithm satisfying \eqref{eq.FPT} is \textit{fixed-parameter tractable (FPT) with respect to $i$}.

For computations of knot invariants, $i$ is often a graph-theoretic width of planar diagrams of knots. Some relevant examples of widths are \textit{tree-width}, \textit{cut-width}, \textit{carving-width}, and, formulated in this paper, \textit{contraction-width}. Among them, tree-width is the most extensively studied, hence we list its relationship with other widths in \cref{tab:width}; we refer to other sources such as \cite{DPSS:treewidth} for its definition. Contraction-width plays a central role in this paper, and its definition is given in \cref{subsec:width}.

%\begin{definition}
%  \label{def:tree-width}
%  A tree decomposition of a graph $\Gamma = (V,E)$ is a tree $T$ whose vertices consist of subsets $X_i$ of $V$ such that 
%  \begin{enumerate}[(a)]
%    \item $\bigcup_i X_i = V$.
%    \item If $v\in X_i\cap X_j$, then any $X_k$ in the path in $T$ connecting $X_i$ and $X_j$ also contains $v$.
%    \item For every edge $\{v,w\}\in E$, there exists $X_i$ such that $v, w \in X_i$. 
%  \end{enumerate}
%  The width of a tree decomposition is $\max_i |X_i|-1$. The \textit{tree-width} of $\Gamma$ is the minimum width among all possible tree decompositions of $\Gamma$.
%\end{definition}

\begin{table}[ht!]
  \centering
  \setlength{\tabcolsep}{4pt}
  \begin{tabular}{@{}cccc@{}}
            \toprule
            Width & Relation with tree-width & Asymptotics & Reference \\
             \midrule
              Tree-width $tw(\Gamma)$ & itself  & $O(\sqrt{n})$ & \cite{LT:separator}\\
            \midrule
            Cut-width $c(\Gamma)$ & $c(\Gamma) = O(tw(\Gamma)\cdot\Delta(\Gamma)\cdot\log n) $ & $O(\sqrt{n}\log n)$ & \cite{KS:TPC}\\
            \midrule
            Carving-width $cw(\Gamma)$ & $\frac{2}{3}(tw(\Gamma) + 1) \leq cw(\Gamma) \leq \Delta(\Gamma)\cdot (tw(\Gamma) + 1)$ & $O(\sqrt{n})$ & \cite{Bienstock:tree} \\
            \midrule
            Contraction-width $w(\Gamma)$ & $\frac{2}{3}(tw(\Gamma) + 1) \leq w(\Gamma) \leq \Delta(\Gamma)\cdot (tw(\Gamma) + 1)$ & $O(\sqrt{n})$ & \eqref{eq:w-tw} and \eqref{eq.w-cw}\\
            \bottomrule
  \end{tabular}
  \vspace{10pt}
  \caption{Different widths and their relationship with the tree-width, where $\Gamma$ is a planar graph with $n$ vertices and maximum degree $\Delta(\Gamma)$.}
  \label{tab:width}
\end{table}

For the FPT algorithms discussed in this paper, we have $f(i) = e^{O(i)}$ and $i = O(\sqrt{n})$ in terms of \cref{eq.FPT}. Previously, FPT algorithms were first established for the Jones polynomial by Makowsky \cite{Makowsky:treewidth} and then for all Reshetikhin--Turaev knot polynomials by Maria \cite{Maria:FPT}, the latter being FPT with respect to the carving-width. Since finding the optimal carving decomposition for planar diagrams is polynomial-time \cite{ST:ratcatcher}, the algorithm in \cite{Maria:FPT} gives a sub-exponential bound $e^{O(\sqrt{n})}$ for computing Reshetikhin--Turaev knot polynomials. However, the algorithm in \cite{Maria:FPT} involves classifying morphism compositions into 7 elementary pieces and dealing with them case by case. To the best of the author's knowledge, no publicly available implementation of the algorithm in \cite{Maria:FPT} exists as of this writing. 

On the other hand, the computation of Reshetikhin--Turaev knot polynomials naturally translates into tensor contractions on the corresponding tensor networks. Tensor networks have been intensively studied in recent years \cite{MS:tensor,Orus:tensor}, and the implementation of tensor contractions is straightforward since its computation is locally uniform. Experiments have revealed tensor contractions to be quite efficient for numerical evaluation of the Jones polynomial at specific values \cite{MK:Jones} and for exact computation of knot polynomials \cite{GL:patterns}. In joint work with S. Garoufalidis, the author gave an implementation in \cite{GL:patterns}, using which the $V_2$-polynomial (also known as the $2$-colored Links--Gould polynomial \cite{Vn}) was computed for millions of knots, including all those up to $16$ crossings; the computed data confirmed that, for all 352 million knots up to $19$ crossings,  the $V_2$-polynomial detects the genus of knots. However, no thorough analysis of the time complexity of the tensor network approach for computing knot polynomials has appeared in the literature, as people in the field of tensor networks usually only consider tensors consisting of scalars instead of polynomials. 

In this paper, we build from the very basics of tensor networks and show that the computation of Reshetikhin--Turaev knot polynomials via tensor contractions is in fact FPT with respect to the contraction-width, which recovers the sub-exponential bound given in \cite{Maria:FPT} with a simpler proof. 

More precisely, let $RT$ be a $k$-variable Reshetikhin--Turaev knot polynomial defined by an $R$-matrix $R\in \Aut(V\otimes V)$ and a ribbon element $h\in \Aut(V)$, and let $d$ be the dimension of $V$. We show the following: 

\begin{theorem}
  \label{th:tractable}
  Given any oriented knot diagram $D$ with $n$ crossings and a contraction sequence $\mathcal{S}$ on $D$ with contraction-width $w$, computing $RT(D)$ by contracting along $\mathcal{S}$ is 
  \begin{equation}
    \label{eq.tractable}
  O(d^{w}n^{k+2}(\log n)^3).
  \end{equation} 
\end{theorem}

\begin{figure}[hbt!]
  \begin{minipage}{0.45\textwidth}
    \centering
    \includegraphics[width=\linewidth]{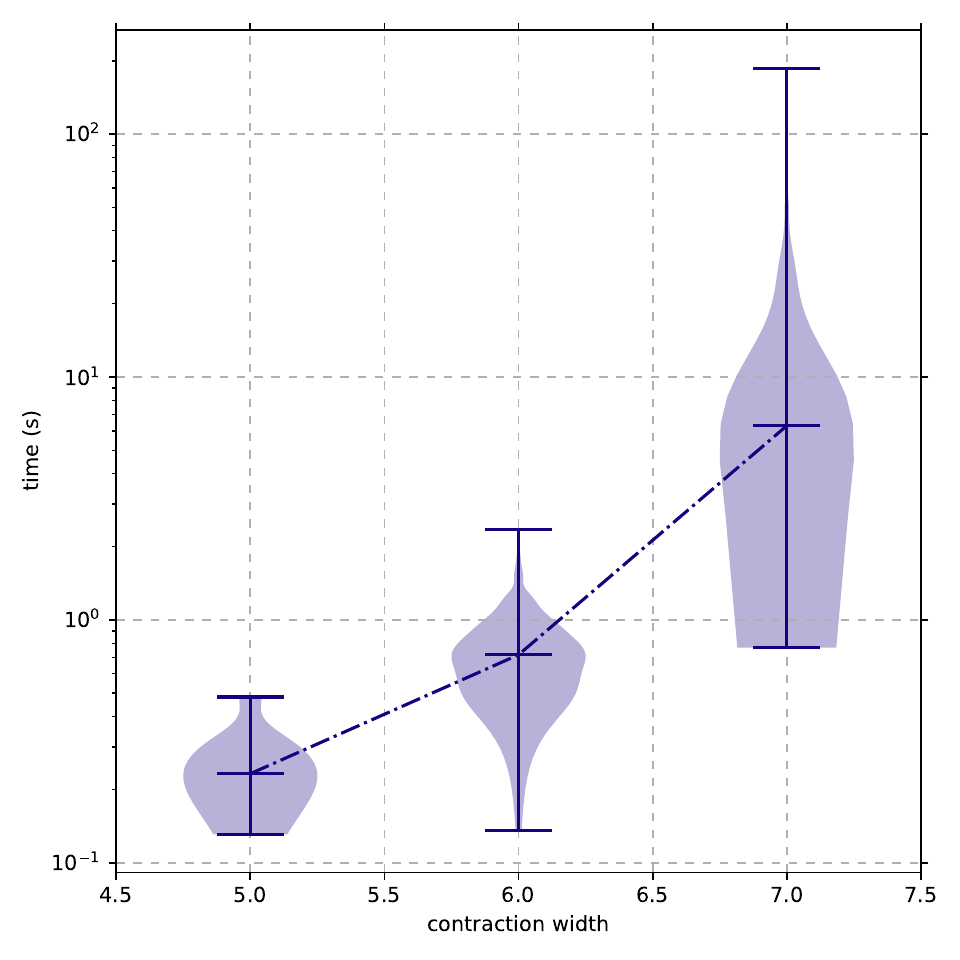}
    \subcaption{The $V_3$-polynomial for the first 1000 knots with 13 crossings in the knot table \cite{HTW}, plotted in log scale.}
    \label{fig:time-13c}
  \end{minipage}
  \hfill
  \begin{minipage}{0.45\textwidth}
    \centering
    \includegraphics[width=\linewidth]{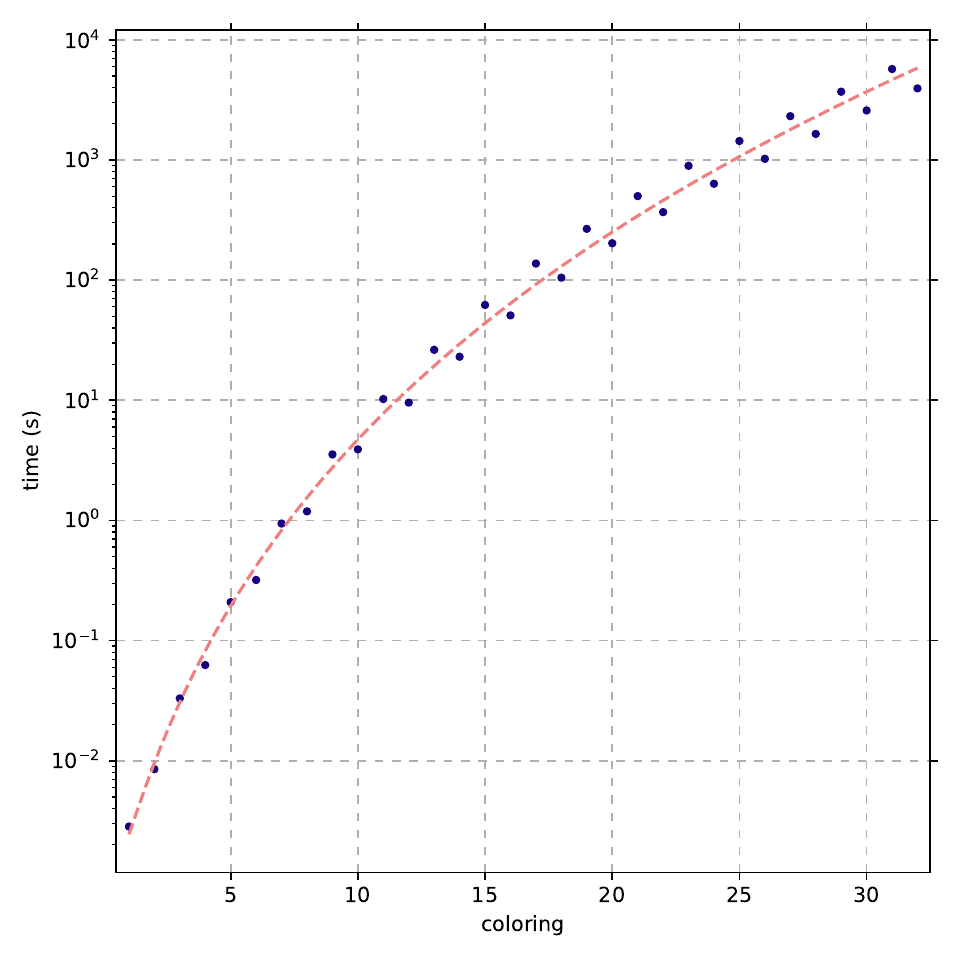}
    \subcaption{The $n$-colored Jones polynomial for the 11n34 knot as $n$ ranges from $1$ to $32$, using a contraction sequence of width $7$, plotted in log scale.}
    \label{fig:time-width6}
  \end{minipage}

    \vspace{1em}  

    \begin{minipage}{0.45\textwidth}
    \centering
    \includegraphics[width=\linewidth]{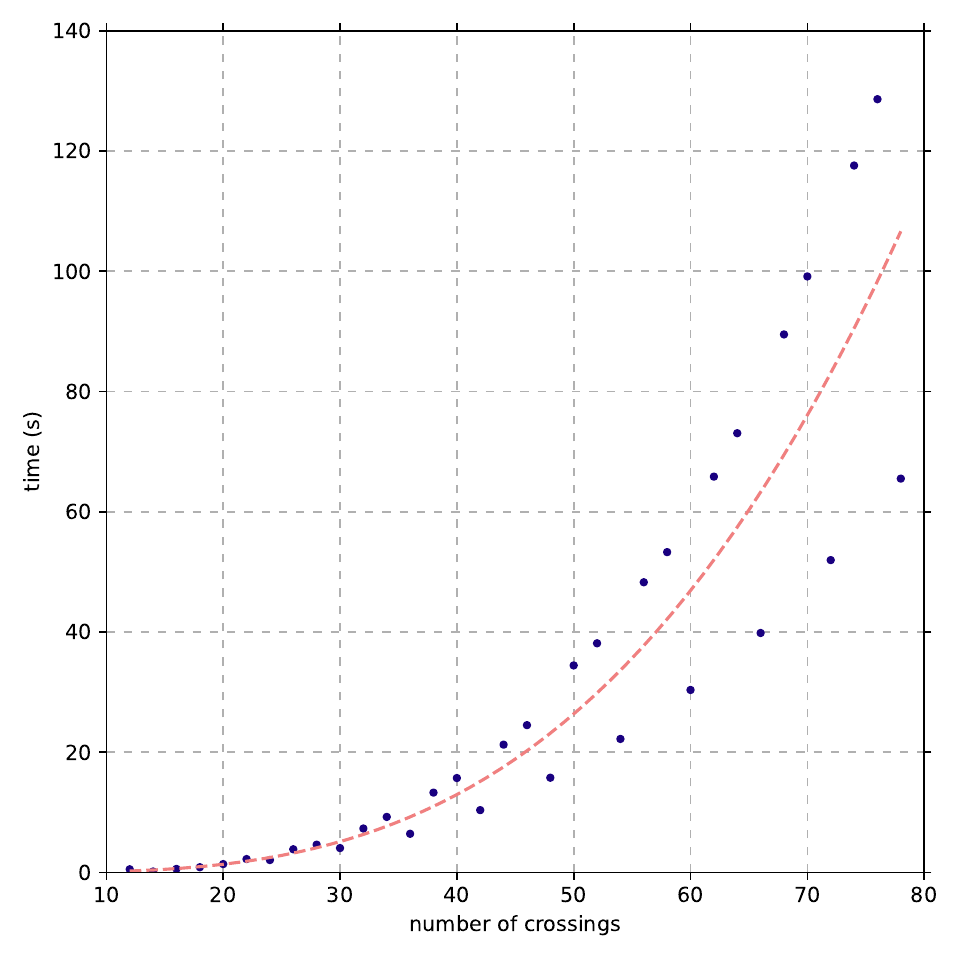}
    \subcaption{The $V_2$-polynomial for $(3,i)$-torus links with $i=6,\dots,39$, all using contraction sequences of width $7$.}
    \label{fig:time-torus}
  \end{minipage}
  \hfill
  \begin{minipage}{0.45\textwidth}
    \centering
    \includegraphics[width=\linewidth]{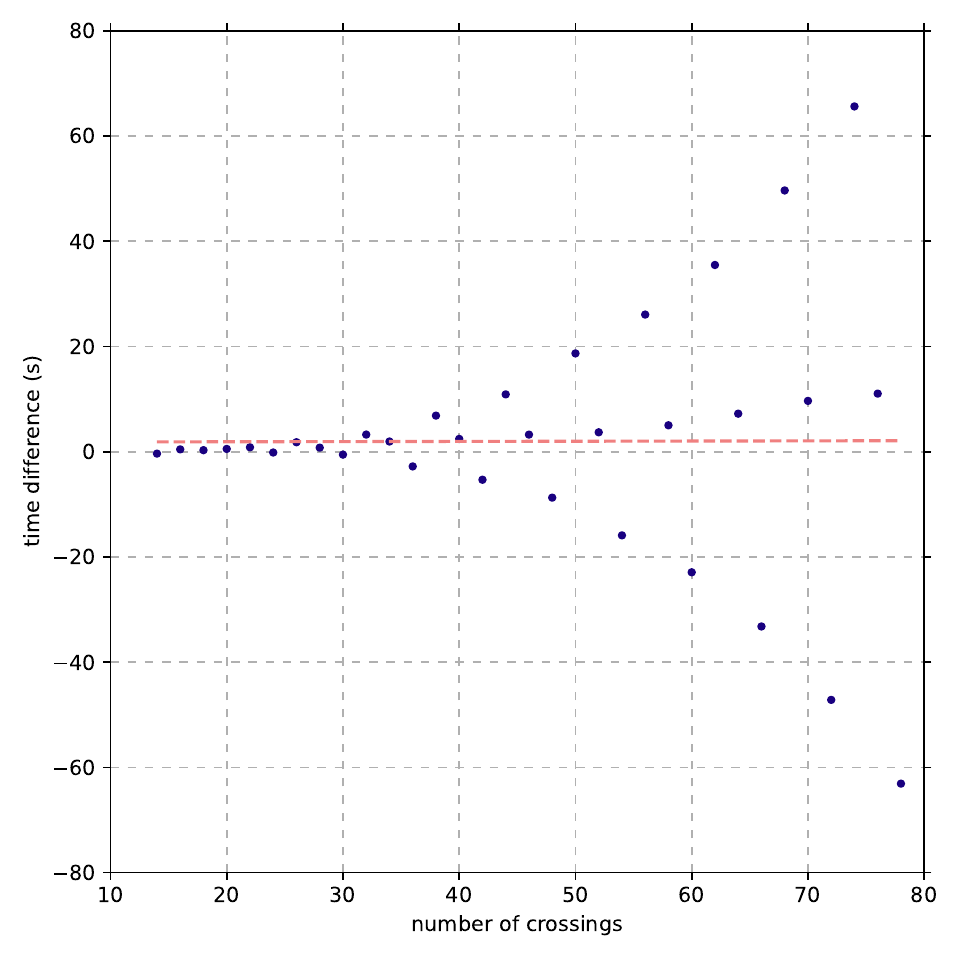}
    \subcaption{Consecutive differences of the computing times in \cref{fig:time-torus}.}
    \label{fig:time-torus-diff}
  \end{minipage}
  \caption{Run time for computing polynomials on various sets of knots. For the $V_n$-polynomial, $k =2$ and $d=4n$; for the $n$-colored Jones polynomial, $k=1$ and $d =n+1$.}
  \label{fig:time-plot}
\end{figure}

We give in \texttt{SnapPy} \cite{snappy} the first general-purpose implementation of this computation, which is now available in the GitHub-version of \texttt{Spherogram} \cite{spherogram}, a submodule of \texttt{SnapPy}. 
%\footnote{The author acknowledges the usage of Claude Code Opus 5.5 in assist of part of the implementation; the implementation is now available in the GitHub-version of \texttt{Spherogram}\cite{spherogram}, a submodule of \texttt{SnapPy}.}.
We benchmarked the run time of our implementation on various sets of knots and knot polynomials, producing \cref{fig:time-plot}. To highlight the capability of our implementation, we also present \cref{table:censusknots,table:randomknots} showing its performance on computing the Jones polynomial on knot diagrams with hundreds and thousands of crossings; see \cref{sec:benchmark} for detailed explanations.

\begin{table}[!htbp]
  \centering
    \begin{tabular}{ccccc}
    \toprule
    Name  & Crossings & Width & Time (minutes) & Memory (GB)\\
    \midrule
    K9\_24 & 213 & 27 & 41.20 & 0.46\\
    K9\_25 & 245 & 28 & 89.85 & 0.50\\
    K9\_27 & 237 & 28 & 52.10 & 0.63\\ 
    K9\_32 & 244 & 25 & 26.69 & 0.40\\
    K9\_33 & 260 & 25 & 24.59 & 0.46\\
    K9\_41 & 266 & 31 & 434.84 & 1.30\\
    K9\_43 & 284 & 32 & 565.27 & 1.45\\
    K9\_44 & 262 & 30 & 160.44 & 1.65\\
    K9\_45 & 288 & 31 & 462.35 & 2.42\\
    \bottomrule
    \end{tabular}
  \vspace{10pt}
  \caption{Run time and RAM consumption of our \texttt{SnapPy} implementation on computing the Jones polynomial of the census knots in \cref{eq:censusknots}, where the width denotes the contraction-width of the contraction sequence being used. In comparison, \texttt{Regina} went OOM for every one of them under a memory constraint of 160GB.}
  \label{table:censusknots}
\end{table}

\begin{table}[!htbp]
  \centering
    \begin{tabular}{cccccc}
    \toprule
    Knot & $O(n^3)$-width & $O(n^2)$-width & \texttt{Regina} & $O(n^3)$-greedy & $O(n^2)$-greedy\\
    \midrule
    rd\_510 & 14 & 15 & (71.74, 16.70) &  (6.05, 0.30) & (0.63, 0.30)\\
    rd\_520 & 14 & 14 &(0.67, 0.55) & (6.17, 0.30) & (0.53, 0.30) \\
    rd\_530 & 16 & 16 &(13.57, 3.67) & (7.36, 0.33) & (1.02, 0.33)\\
    rd\_540 & 18 & 18 &(263.68, 75.94) & (8.10, 0.36) & (1.97, 0.36)\\
    rd\_550 & 15 & 17 &(14.25, 6.55) & (7.38, 0.31) & (1.18, 0.32)\\
    \dashedmidrule{6}
    rd\_950 & 22 & 24 &OOM & (45.41, 0.77) & (106.98, 2.49)\\
    rd\_960 & 15 & 15 &(18.80, 1.90) & (33.90, 0.32) & (1.74, 0.32)\\
    rd\_970 & 19 & 20 &OOM & (42.02, 0.47) & (7.11, 0.41)\\
    rd\_980 & 20 & 20 &OOM & (49.68, 0.68) & (9.97, 0.47)\\
    rd\_990 & 16 & 16 &(14.41, 4.59) & (38.06, 0.36) & (2.70, 0.36)\\
    rd\_1000 & 17 & 17 &(135.97, 24.48) & (41.75, 0.41) & (5.78, 0.40)\\
    \dashedmidrule{6}
    rd\_2010 & 17 & 20 &OOM & (307.89, 0.50) & (6.57, 0.39)\\
    rd\_2012 & 25 & 26 &OOM & (495.01, 6.42) & (649.65, 9.79) \\
    rd\_2025 & 19 & 20 &OOM & (272.21, 0.70) & (12.51, 0.69)\\
    rd\_3037 & 21 & 21 &OOM & (996.90, 4.49) & (119.06, 2.77)\\
    rd\_3070 & 23 & 26 &OOM & (1316.62, 4.23) & (1024.28, 10.71)\\
    rd\_4012 & 26 & 27 &OOM & (5221.20, 76.78) & (3813.91, 55.74) \\
    %rd\_5005 & & &OOM & \\
    %rd\_6037 & & &OOM & \\
    %rd\_8061 & & & \\
    %rd\_10147 & & OOM & \\
    %rd\_17358 & & OOM & \\
    \bottomrule
    \end{tabular}
  \vspace{10pt}
  \caption{Pairs of (run time (minutes), RAM consumption (GB)) on computing the Jones polynomial within 6 days and 160GB RAM, where $\mathrm{rd}\_n$ denotes the generic random knot with $n$ crossings in our collection and the last two columns list the pairs for our \texttt{SnapPy} implementation. See \cref{subsec:jones} for explanation of the $O(n^3)$- and $O(n^2)$-greedy algorithms.}
  \label{table:randomknots}
\end{table}

We remark that the generality and efficiency of our implementation enable investigations to much broader families of knot polynomials than people could have done before, for example the $\mathrm{FK}3$ polynomial that S. Garoufalidis and the author defined and computed in \cite{GL:FK3}.

A bit more discussion is needed to recover the sub-exponential bound in \cite{Maria:FPT}. We will need
the following consequence of a result by Markov and Shi \cite{MS:tensor}, explained in \cref{subsec:width}:

\begin{lemma}
  \label{le:bd-tw}
  For any graph $\Gamma$, 
  \begin{equation}
    \label{eq:w-tw}
    w(\Gamma) \leq \Delta(\Gamma)\cdot(tw(\Gamma) + 1),
  \end{equation}
  where $\Delta(\Gamma)$ is the maximum degree of $\Gamma$.
\end{lemma}

Applying an approximation algorithm for generating contraction sequences (see \cref{subsec:co-proof}), we can drop the requirement of the contraction sequence from \cref{th:tractable}, obtaining: 

\begin{corollary}
  \label{co:tractable}
  For any oriented knot diagram $D$ with $n$ crossings, computing $RT(D)$ is 
  \begin{equation}
  O(d^{\frac{3}{2}w(D)}n^{k+2}(\log n)^3).
  \end{equation}
\end{corollary}

Applying \cref{le:bd-tw} and substituting $tw(D)$ with $O(\sqrt{n})$, \cref{co:tractable} translates into
\begin{corollary}
  \label{co:sqrtn}
  %Given a knot polynomial $RT$ of Reshetikhin--Turaev type. 
  For any oriented knot diagram $D$ with $n$ crossings, $RT(D)$ can be computed in time $e^{O(\sqrt{n})}$.
\end{corollary}
This recovers the sub-exponential bound given in \cite{Maria:FPT}. 

In fact, \cref{th:tractable} is an improvement to the results in \cite{Maria:FPT} with a different approach resulting in a simpler proof. This is because of the following relationship between the contraction-width $w(\Gamma)$ and carving-width $cw(\Gamma)$:
\begin{equation}
  \label{eq.w-cw}
  cw(\Gamma) \leq w(\Gamma) \leq \frac{3}{2} cw(\Gamma),
\end{equation}
where the inequalities may be strict. The above inequality is due to Bienstock \cite[Remark 3]{Bienstock:tree}, in which vertex congestion coincides with our contraction-width and it is well-known that edge congestion is equal to the carving-width. The proof of \cref{eq.w-cw} is elementary using the definitions of vertex and edge congestions;
we refer the reader to the work of O'Gorman \cite{OGorman:FPT} for how vertex and edge congestions relate to the costs of tensor network contraction, which, along with our discussions here, would demonstrate how the vertex congestion coincides with our contraction-width.

We give a brief review of the Reshetikhin--Turaev functor in \cref{sec:setup}, explaining how the computation translates into tensor contractions, and then give the definition of contraction-width and detailed analysis for the time complexity in \cref{sec:time}, proving \cref{th:tractable} along with other results above. We end this paper with benchmarks comparing our implementation with existing publicly available packages in \cref{sec:benchmark}, showing that our implementation is several orders of magnitude more efficient than others. Throughout this paper, Reshetikhin--Turaev knot polynomials are assumed to be given by $R$-matrices and ribbon elements with matrix representations whose entries are Laurent polynomials with integer coefficients.

We remark that our results easily generalize to links and possibly with differently colored components\footnote{Our implementation currently supports links with uniformly colored components.}, but we talk about knots exclusively for simplicity and because ``knot polynomial'' is a much more common phrase than ``link polynomial''.

\section{Setup}
\label{sec:setup}

We give some basic descriptions about our settings in this section, laying the groundwork for the discussion of time complexity in \cref{sec:time}. 

\subsection{Long knot diagrams and rotation numbers}
\label{subsec:longknot}

Given a planar diagram of an oriented knot, we convert it into an oriented long knot diagram by cutting an arc open and pulling the two open ends towards infinity in opposite directions without creating any self-intersection. Following the convention of Bar-Natan and van der Veen in \cite{BNV:API}, we regard the out-pointing strand to be pulled vertically upward, and the other one vertically downward. We further require that the arcs always point upward in small neighborhoods of each crossing by applying an ambient isotopy to rotate the crossings. 

%Figure~\ref{f.long4-1} presents an oriented long knot diagram of the $4_1$ knot. 

\begin{figure}[htpb!]
\centering
\scalebox{.6}{$\vcenter{\hbox{\includegraphics{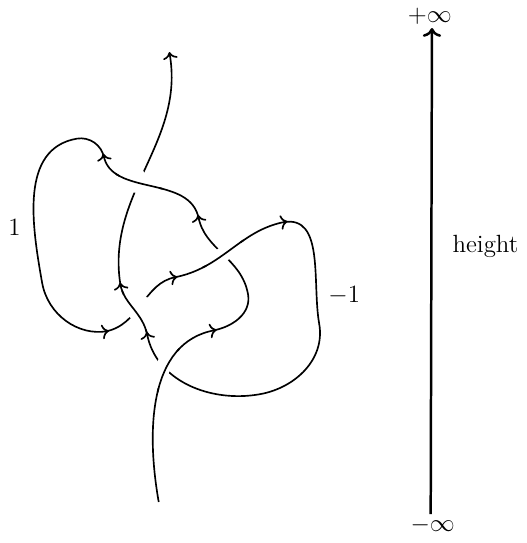}}}$
}
\caption{An oriented long knot diagram of the $4_1$ knot, with nonzero
  rotation numbers labeled.}
\label{f.long4-1}
\end{figure}

The vertical axis now gives us a natural height function, which
restricts to a Morse function along the long knot diagram, whose local maxima
and minima occur only on the arcs and never at the crossings. The rotation number associated to each arc in the oriented long knot diagram is defined as follows.

\begin{definition}
The \textit{rotation number} of an arc $\mathcal{A}$ in an oriented long knot
diagram equipped with a height function $\ell$ as described above, is the integer
\begin{equation*}
  \sum_{p\in \{d\ell = 0 \} \cap \mathcal{A}} (-1)^{\delta_\ell(p)} \varepsilon_\mathcal{A} (p),
\end{equation*}
where $\delta_\ell(p) = 1$ if $p$ is a local maximum of $\ell$ and $0$ otherwise;
$\varepsilon_\mathcal{A} (p) = 1$ if the oriented arc $\mathcal{A}$ points to the
right at $p$, and $0$ otherwise. 
\end{definition}

Note that the rotation number is determined not only by the combinatorial information of the oriented long knot diagram, but also by how the diagram is embedded in the plane, inducing the height function. Given a planar diagram of an oriented knot with $n$ crossings, there is an $O(n)$ algorithm that converts it into an oriented long knot diagram with rotation numbers in $\{0,\pm 1\}$ \cite{BNV:API}. A description of the algorithm can be found in \cite[Section 2.2.1]{GL:patterns}.

\subsection{The Reshetikhin--Turaev functor}
\label{subsec:RT}

Depending on the context, the precise constructions of the Reshetikhin--Turaev functor vary, but they all share some common structures, of which only those relevant to computation are described here; we refer the reader to \cite{Ohtsuki:quantum,GK:multi} for introductions to the theoretical aspects.

The central idea of the Reshetikhin--Turaev functor is to view oriented long knot diagrams as an analogue of $1$-dimensional cobordisms, assembled from elementary pieces, and to map the pieces into the category of vector spaces. The composition of the linear maps obtained by applying the functor is thus the desired knot invariant. 

Therefore, as with any topological quantum field theory, the Reshetikhin--Turaev functor is determined by its image on those elementary pieces. The information needed consists of:
\begin{itemize}
  \item A vector space $V$ over a field $\mathbb{F}$.
  \item An \textit{$R$-matrix} $R\in \Aut(V\otimes V)$.
  \item A \textit{ribbon element} $h \in \Aut(V)$.
\end{itemize}
We remark that in practice $R$ is sparse due to weight conservation \cite{Ohtsuki:quantum}, which our implementation heavily exploits.
The Reshetikhin--Turaev functor is then determined as follows:

\begin{equation*}
  \begin{aligned}
    &\vcenter{\hbox{\includegraphics{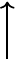}}} 
    \mapsto
    \begin{tikzcd}
    V \\
    V\arrow[u,swap, "1_V"]
    \end{tikzcd}
    \quad
   &\vcenter{\hbox{\includegraphics{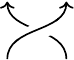}}} 
    \mapsto 
    \begin{tikzcd}
    V\otimes V \\
    V\otimes V\arrow[u,swap, "R"]
    \end{tikzcd}
    \quad
    &\vcenter{\hbox{\includegraphics{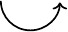}}}
    \mapsto
    \begin{tikzcd}
    V^*\otimes V\\
    \mathbb{F} \arrow[u,swap, "u"]
    \end{tikzcd}
    \quad
    &\vcenter{\hbox{\includegraphics{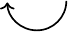}}}
    \mapsto
    \begin{tikzcd}
    V\otimes V^*\\
    \mathbb{F}\arrow[u,swap, "u'"]
    \end{tikzcd}
    \\
    &\vcenter{\hbox{\includegraphics{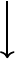}}} 
    \mapsto
    \begin{tikzcd}
    V^* \\
    V^*\arrow[u,swap, "1_{V^*}"]
    \end{tikzcd}
    \quad
    &\vcenter{\hbox{\includegraphics{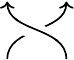}}}
    \mapsto
    \begin{tikzcd}
    V\otimes V \\
    V\otimes V\arrow[u,swap, "R^{-1}"]
    \end{tikzcd}
    \quad
    &\vcenter{\hbox{\includegraphics{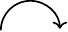}}}
    \mapsto
    \begin{tikzcd}
    \mathbb{F} \\
    V\otimes V^*\arrow[u,swap, "n"]
    \end{tikzcd}
    \quad
    &\vcenter{\hbox{\includegraphics{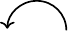}}}
    \mapsto
    \begin{tikzcd}
    \mathbb{F} \\
    V^*\otimes V\arrow[u,swap, "n'"]
    \end{tikzcd}
  \end{aligned}
\end{equation*}
where $u$ and $u'$ are defined by
\begin{equation*}
  u(1) = \sum_i e^*_i\otimes h(e_i),\quad u'(1) = \sum_i e_i\otimes e_i^*,
\end{equation*}
where $\{e_i\}$ is any basis of $V$ and $\{e_i^*\}$ is the dual basis of $\{e_i\}$, and $n$ and $n'$ are defined by
\begin{equation*}
  n(x\otimes f) = f(h^{-1}(x)), \quad n'(f\otimes x) = f(x)
\end{equation*}
for any $x\in V$ and $f \in V^*$.

Since $(n\otimes 1_V)\circ (1_V\otimes u) = 1_V$ and $(n'\otimes 1_{V^*})\circ (1_{V^*} \otimes u')  = 1_{V^*}$, it is easy to see that the result of applying the Reshetikhin--Turaev functor is determined by the combinatorial information of the planar diagram along with the rotation numbers. With the rotation numbers living in $\{0,\pm1\}$, the result can be described by $O(n)$ elementary pieces. Hence it is $O(n)$ to convert an oriented long knot diagram into compositions of linear maps as above. 

The result of applying the Reshetikhin--Turaev functor to an oriented long knot diagram is a diagonalizable linear map in $\End(V)$ whose eigenvalues are all the same. To obtain a knot polynomial, indeterminates are plugged into the definitions of $R$ and $h$, turning the resulting linear map into a parametrized family, and the eigenvalue into the polynomial. 

\subsection{Tensor networks and contractions}
\label{subsec:tensor}

We briefly review some standard terminology of tensor networks and tensor contractions. We refer the reader to \cite{GL:patterns} for slightly more detailed descriptions, and to \cite{JC:tensornet} for an extensive and graphical introduction to this topic.

\begin{definition}
  Let $\mathcal{R}$ be a (commutative) ring. An \textit{$n$-tensor} $T$ over $\mathcal{R}$ is a tuple
  $(T_{i_1,\dots,i_n})_{(i_1,\dots,i_n)\in \mathcal{S}}$ where $T_{i_1,\dots,i_n}\in \mathcal{R}$
  and the index set $\mathcal{S}$ is of the form 
\begin{equation*}
  \mathcal{S} = \left\{ 1,\dots,m_1 \right\}\times \dots\times \left\{
    1,\dots, m_n \right\} \subset \mathbb{N}^n,
\end{equation*}
We call the integer $m_k$ ($k\in \left\{ 1,\dots,n \right\}$) the dimension
of the \textit{leg} $k$ of the tensor $T$. 
\end{definition}

\begin{definition}
A \textit{tensor network} over a ring $\mathcal{R}$ is a (multi)graph $\Gamma$ where each vertex $v$ is associated with
a tensor $T_v$ over $\mathcal{R}$ whose legs are bijectively associated with the ends of the edges incident
to $v$. Such an association must satisfy that for each edge in $\Gamma$, its two ends are associated with two different
legs with the same dimension of one or two tensors.
\end{definition}

Given a tensor network $\Gamma$ over $\mathcal{R}$, for an $n$-tensor $T$ and an $n'$-tensor $T'$ in $\Gamma$,
if two legs $k$ of $T$ and $k'$ of $T'$ are associated with the two ends of an edge $e$ in $\Gamma$, we can
\textit{contract $T$ and $T'$ along $e$} to obtain an $(n+n'-2)$-tensor $T''$, defined by
\begin{equation*}
T''_{i_1,\dots,\widehat{i_k},\dots,i_n,j_1,\dots,\widehat{j_{k'}},\dots,j_{n'}}
\coloneqq \sum_{\substack{i_k = j_{k'} \in \left\{ 1,\dots,m_k \right\}
%\\ i_k = j_{k'} 
}} T_{i_1,\dots,i_k,\dots,i_n}T'_{j_1,\dots,j_{k'},\dots,j_{n'}}
\end{equation*}
where the hats indicate that the corresponding indices are deleted. Moreover, we can remove the edge $e$ and replace $T$ and $T'$ and their corresponding vertices with $T''$ and a new vertex, obtaining a tensor network with one fewer vertex and edge than $\Gamma$. 

More generally,
if legs $k_1,\dots,k_s$ of $T$ and legs
$k'_1,\dots,k'_s$ of $T'$ are associated with the two ends of edges $e_1,\dots,e_s$ in $\Gamma$ respectively, we can contract $T$ and $T'$ along $e_1,\dots,e_s$ at once to obtain an $(n+n' - 2s)$-tensor $T''$ defined by
\begin{equation}
\label{tensorcontract}
T''_{i_1,\dots,\widehat{i_{k_1}},\dots,
  \widehat{i_{k_s}},\dots,i_n,j_1,\dots,\widehat{j_{k'_1}},\dots,
  \widehat{j_{k'_s}},\dots,j_{n'}}
\coloneqq \sum_{ 
\substack{i_{k_1} = j_{k'_1} \in \{ 1,\dots,m_{k_1}\}
\\ \vdots\\
i_{k_s} = j_{k'_s} \in \{ 1,\dots,m_{k_s} \}
}
} T_{i_1,\dots,i_n}T'_{j_1,\dots,j_{n'}}
\end{equation}
We can then remove the edges $e_1,\dots,e_s$ and replace $T$ and $T'$ and their corresponding vertices with $T''$ and a new vertex, obtaining a tensor network with one fewer vertex and $s$ fewer edges than $\Gamma$. Similarly, we can define contractions of more than two tensors at once.

\begin{example}
  \label{ex.RTnet}
  Given an $R$-matrix $R\in \Aut(V\otimes V)$ and a ribbon element $h\in \Aut(V)$ determining a Reshetikhin--Turaev functor, we fix a basis $\mathcal{B}\coloneqq \{e_1,\dots,e_{d}\}$ of $V$. Under $\mathcal{B}$, $R^{\pm1}$ becomes a $4$-tensor $(R^{\pm 1})_{e_i\otimes e_j}^{e_k\otimes e_l}$ as $(i,j,k,l)$ varies in $\{1,\dots,d\}^4$. Similarly $h^{\pm1}$ becomes a $2$-tensor $(h^{\pm 1})_{e_i}^{e_j}$. For an oriented long knot diagram with $n$ crossings, the Reshetikhin--Turaev functor turns it into a tensor network by the following assignment:
\begin{equation*}
  \begin{aligned}
    &\vcenter{\hbox{\includegraphics{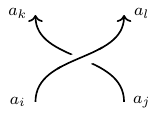}}}  \mapsto 
    \vcenter{\hbox{\includegraphics{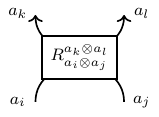}}} 
    \quad
    &\vcenter{\hbox{\includegraphics{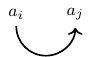}}}
    \mapsto
    \vcenter{\hbox{\includegraphics{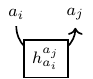}}}\\
    &\vcenter{\hbox{\includegraphics{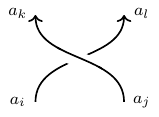}}}
    \mapsto  
    \vcenter{\hbox{\includegraphics{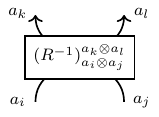}}}
    &\vcenter{\hbox{\includegraphics{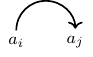}}}
    \mapsto
    \vcenter{\hbox{\includegraphics{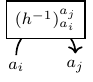}}}
  \end{aligned}
\end{equation*}
That is, it turns each crossing into a degree $4$ vertex associated with a $4$-tensor $(R^{\pm 1})_{e_i\otimes e_j}^{e_k\otimes e_l}$, and each cup or cap going rightwards into a degree $2$ vertex associated with a $2$-tensor $(h^{\pm 1})_{e_i}^{e_j}$. With rotation numbers in $\{0,\pm 1\}$, we obtain a tensor network with no more than $2n$ vertices, where all legs are of dimension $d$.
\end{example}

\begin{figure}[htpb!]
\centering
    \begin{equation*}
\vcenter{\hbox{\includegraphics{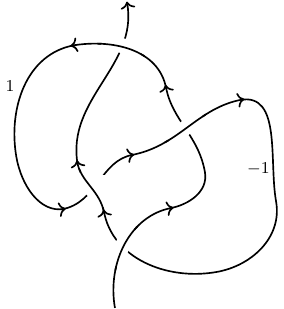}}}
    \quad\rightsquigarrow \quad 
    \vcenter{\hbox{\includegraphics{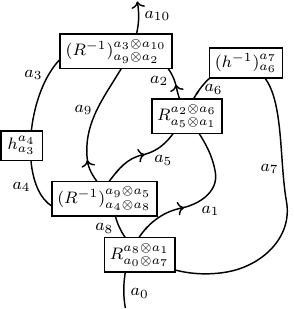}}}
\end{equation*}
\caption{The tensor network obtained from an oriented long knot diagram
      of the $4_1$ knot.}
\label{f.tensornet4-1}
\end{figure}

\begin{example}
  \label{ex.4-1-v}
  \cref{f.tensornet4-1} illustrates the tensor network obtained by applying the Reshetikhin--Turaev functor on an oriented long knot diagram of the $4_1$ knot.
  With the symbols in \cref{f.tensornet4-1}, one can check that the corresponding knot polynomial, as discussed at the end of \cref{subsec:RT}, can be expressed as the following finite sum:
  \begin{equation}
    \label{eq.4-1-RT}
    \sum_{\substack{a_1,\dots,a_9\in \mathcal{B}\\ a_0 = a_{10} = e_1}}
    R_{a_0\otimes a_7}^{a_8\otimes a_1} \cdot 
    h_{a_3}^{a_4}\cdot  
    (R^{-1})_{a_4\otimes a_8}^{a_9\otimes a_5} \cdot 
    R_{a_5\otimes a_1}^{a_2\otimes a_6} \cdot
    (h^{-1})_{a_6}^{a_7} \cdot
    (R^{-1})_{a_9\otimes a_2}^{a_3\otimes a_{10}}.
  \end{equation}
  Recall that $\mathcal{B} = \{e_1,\dots,e_{d}\}$ is the basis of $V$.
\end{example}

\cref{ex.RTnet,ex.4-1-v} demonstrate that all Reshetikhin--Turaev knot polynomials can be expressed as a sum similar to \eqref{eq.4-1-RT}. Since the sum is finite, we can compute it piece by piece in an arbitrary order, which naturally translates into doing a sequence of tensor contractions on the tensor network until only one vertex remains. Such a sequence is called a \textit{contraction sequence}. A contraction sequence can be encoded by an ordered list of sets of edges removed at each step, and it plays a crucial role in the time complexity of computing the contractions, which we discuss in the following section. Note that we can always eliminate two free legs of the tensor network by, for example, setting at the beginning $a_0=a_{10}=e_1$ in \cref{ex.4-1-v}.

\section{Bounding the time complexity}
\label{sec:time}

\subsection{Contraction-width}
\label{subsec:width}

We now restrict ourselves to tensor networks obtained by applying the Reshetikhin--Turaev functor to oriented long knot diagrams, where all legs of tensors have dimension $d = \dim V$.

It is clear from \eqref{tensorcontract} that it takes $O(c\cdot d^s)$ ring operations (i.e. additions and multiplications in $\mathcal{R}$) to compute a single
entry in a tensor resulting from contracting $c$ tensors along $s$ edges at once. Let $l_1,\dots,l_c$ be the number of legs of the $c$ tensors, the resulting tensor thus has $- 2s +\sum_{i=1}^c l_i $ legs, hence consists of $d^{- 2s + \sum_{i=1}^c l_i}$ entries. Therefore the total number of ring operations required to build the resulting tensor is
\begin{equation*}
  O(c\cdot d^s)\cdot d^{- 2s + \sum_{i=1}^c l_i} = O(c\cdot d^{-s + \sum_{i=1}^c l_i}).
\end{equation*}
This motivates the following definition.

\begin{definition}
Given a contraction sequence $\mathcal{S}$ on a tensor network $\Gamma$, for a step in $\mathcal{S}$ that contracts $c$ tensors with $l_1,\dots,l_c$ legs respectively along $s$ edges at once, its \textit{local contraction-width} is $-s + \sum_{i=1}^c l_i$. The \textit{contraction-width} $w(\mathcal{S})$ of $\mathcal{S}$ is the maximum local contraction-width of all its steps, and the contraction-width $w(\Gamma)$ of $\Gamma$ is the minimal contraction-width among all possible contraction sequences on $\Gamma$. 
\end{definition}

Note that the contraction-width is independent of the dimension $d$, hence it makes sense to consider (abstract) contraction sequences on a graph $\Gamma$ and define its contraction-width $w(\Gamma)$ similarly. Note also that adding degree $2$ vertices does not affect the contraction-width of a nontrivial graph. These justify the use of $w(D)$ in \cref{sec.intro}.

Since we would like the local contraction-width to be as small as possible, we would always want to maximize the number of legs contracted within each contraction; in doing so, we see that in terms of the underlying graph of $\Gamma$, the quantity $-s + \sum_{i=1}^c l_i$ is equal to the number of distinct edges associated with the legs of the $c$ tensors involved in the contraction.
It is then clear that adding more tensors to a step of contraction never decreases its local contraction-width, hence at best has no effect on the exponential factor of the number of ring operations required. Therefore we further restrict ourselves to contraction sequences where each step contracts exactly $2$ tensors. In particular, for a tensor network with $2n$ vertices, all contraction sequences consist of $2n-1$ steps. For a contraction sequence $\mathcal{S}$ with width $w$, it then takes $O(n\cdot d^w)$ ring operations in $\mathcal{R}$ to perform all steps in $\mathcal{S}$.

In \cite{MS:tensor}, Markov and Shi restricted themselves to contraction sequences where each step contracts only one edge. For such a contraction sequence $\mathcal{S}$, they defined its \textit{contraction complexity} $cc(\mathcal{S})$ as the maximum degree of a merged vertex during the contraction process, and hence the contraction complexity $cc(\Gamma)$ of a tensor network $\Gamma$ is the minimum complexity among all possible $\mathcal{S}$ where each step contracts along only one edge. It is immediate from definition that $w(\mathcal{S}) = cc(\mathcal{S}) + 1$, hence 
\begin{equation}
  \label{eq:w-leq-cc}
  w(\Gamma)\leq cc(\Gamma) + 1,
\end{equation}
where the inequality may be strict.
%since the contraction-width is defined over a larger family of contraction sequences. 
Recall that the \textit{line graph} $\Gamma^*$ of a graph $\Gamma$ is defined by $V(\Gamma^*) \coloneqq E(\Gamma)$ and 
\begin{equation*}
  E(\Gamma^*) \coloneqq \{\{e,f\}\subset E(\Gamma) \mid e\neq f, \ \exists v\in V(\Gamma) \text{ such that } v\in e, v\in f \}.
\end{equation*}
The following were proved in \cite{MS:tensor}:

\begin{proposition}[Markov--Shi]
  \label{prop:MS-line-graph}
  For any graph $\Gamma$, $cc(\Gamma) = tw(\Gamma^*)$. Furthermore, given a tree decomposition of $\Gamma^*$ of width $w$, there is a polynomial-time algorithm translating it into a contraction sequence $\mathcal{S}$ on $\Gamma$ with $cc(\mathcal{S}) \leq w$.
\end{proposition}

\begin{lemma}[Markov--Shi]
  \label{lem:MS-tw}
  For any graph $\Gamma$ of maximum degree $\Delta(\Gamma)$, 
  \begin{equation}
     (tw(\Gamma) -1)/2 \leq tw(\Gamma^*) \leq \Delta(\Gamma)(tw(\Gamma) + 1)- 1.
  \end{equation}
\end{lemma}
Therefore combining \eqref{eq:w-leq-cc}, \cref{prop:MS-line-graph} and \cref{lem:MS-tw}, we have established that
\begin{equation}
  w(\Gamma) \leq cc(\Gamma) + 1 = tw(\Gamma^*) + 1 \leq \Delta(\Gamma)(tw(\Gamma) + 1).
\end{equation}
This gives \cref{le:bd-tw} as we stated.

We end this section with a remark that the contraction-width is also linearly bounded by the cut-width from above; more precisely, we have $w(\Gamma)\leq 3 c(\Gamma)$, recalling that $c(\Gamma)$ denotes the cut-width of $\Gamma$. This is because that a diagram realizing the graph's cut-width gives a sequential order for the tensors, where the tensors are connected one by one and the legs of each tensor can be partitioned into two sides, where the number of legs on each side does not exceed $c(\Gamma)$. In particular, by fixing a small braid index, one can construct knots with arbitrarily many crossings while maintaining a small contraction-width.

\subsection{Proof of Theorem~\ref{th:tractable}
% and Theorem~\ref{th:tractable-precise}
}
\label{subsec:proof}

Given an oriented knot diagram $D$ with $n$ crossings and a contraction sequence on $D$ with contraction-width $w$, we cut an arbitrary arc in $D$ open, turning it into an oriented long knot diagram $D'$. As discussed in \cref{subsec:RT}, the arc cut open can be eliminated immediately after we apply the Reshetikhin--Turaev functor to $D'$, obtaining a tensor network with at most $2n$ vertices along with a contraction sequence whose width is at most $w$. This procedure is $O(n)$ and is dominated by subsequent computation. 

We have established in the previous section that it takes $O(n\cdot d^w)$ ring operations to compute the desired knot polynomial, hence it remains only to bound the time complexity of performing ring operations in $\mathcal{R}$. 

For a $k$-variable Reshetikhin--Turaev knot polynomial, we have $\mathcal{R} = \mathbb{Z}[t_1^{\pm 1}, \dots, t_k^{\pm1}]$. Let $M-1$ be the maximum degree (i.e. the difference between the highest and lowest exponents) of all variables in the entries of $(R_{e_i\otimes e_j}^{e_k\otimes e_l})$ and $(h_{e_i}^{e_j})$. The maximum degree for each variable that can appear throughout the process is thus $2nM$. Hence all polynomials have at most $(2nM)^k = O(n^{k})$ terms. Therefore the naive number of ring operations in $\mathbb{Z}$ required for a single polynomial addition or multiplication is $O(n^{2k})$. Alternatively, we may use Kronecker substitution to turn multivariable polynomials into univariate ones and apply fast Fourier transform to compute their multiplications, which reduces the required number of ring operations in $\mathbb{Z}$ to $O(k(4nM+1)^k\log (4nM+1)) = O(n^k\log n)$ \cite[Section 8.4]{ModernComputerAlgebra}.

We now bound the integers that may appear throughout the process. Let $C$ be the maximum absolute value of coefficients of the entries of $(R_{e_i\otimes e_j}^{e_k\otimes e_l})$ and $(h_{e_i}^{e_j})$. Since there are at most $M^k$ distinct monomials appearing in the entries and the whole state sum is to sum up $d^{2n}$ terms where each term is a product of $2n$ polynomials, the maximum absolute value of integers that can appear throughout the process is bounded by the maximum absolute value of the coefficients of
\begin{equation}
  d^{2n} \left(\sum_{i=1}^{M^k} C m_i\right)^{2n} = \sum_{\substack{p_1 + \cdots + p_{M^k} = 2n\\ p_1,\dots,p_{M^k}\geq 0}} \frac{d^{2n}\cdot (2n)!\cdot C^{2n}}{(p_1!)\cdots (p_{M^k}!)} \cdot m_1^{p_1}\cdots m_{M^k}^{p_{M^k}},
\end{equation}
where each $m_i$ represents a monic monomial. Hence the maximum absolute value of integers is bounded by
\begin{equation}
  \sum_{\substack{p_1 + \cdots + p_{M^k} = 2n\\ p_1,\dots,p_{M^k}\geq 0}}\frac{d^{2n}\cdot (2n)!\cdot C^{2n}}{(p_1!)\cdots (p_{M^k}!)}\leq  \binom{2n + M^k}{M^k} \cdot (2ndC)^{2n} \leq (2n dC M^k)^{2n+M^k},
\end{equation}
whose bit-size is $O(n\log n)$.
Therefore \cref{th:tractable} follows from the fact that adding and multiplying integers with bit-size in $O(n\log n)$ are
\begin{equation}
  O((n\log n)\cdot \log (n\log n)) = O(n(\log n)^2).
\end{equation}

\cref{th:tractable} states only the bound of the time complexity. We remark here that one sees by putting the above analysis together that the memory size of the largest tensor appearing throughout the procedure is also FPT with respect to $w$. Since we need only to store at most $2n+1$ tensors at a time, the memory bound of our computation is also FPT w.r.t. $w$.

\subsection{Proof of Corollary~\ref{co:tractable}}
\label{subsec:co-proof}

We have established that in order to efficiently compute Reshetikhin--Turaev knot polynomials with tensor networks, one should produce on any given planar diagram $D$ a contraction sequence $\mathcal{S}$ whose width is as close to $w(D)$ as possible. We say that a contraction sequence $\mathcal{S}$ on $D$ is \textit{optimal} if $w(\mathcal{S}) = w(D)$. The task of finding an optimal contraction sequence on a given graph should be $\mathrm{NP}$-hard since it includes the task of computing the optimal vertex congestion, where the latter was claimed to be $\mathrm{NP}$-hard in \cite{Bienstock:tree} without a proof; the hardness restricted to (4-valent) planar graphs remains open to the best of the author's knowledge. 

In practice, optimized algorithms of breadth-first or depth-first searches \cite{RJF:optimaltensor}, albeit scale exponentially, are able to produce optimal contraction sequences for knot diagrams up to $30$ or more crossings; for instance, in our implementation in \texttt{SnapPy}, we use the package \texttt{opt\_einsum} \cite{opteinsum} to find contraction sequences. In theory, to obtain a statement such as \cref{co:tractable}, one can use approximation algorithms for contraction sequences instead of optimal ones, which, for any input graph $\Gamma$, gives a contraction sequence $\mathcal{S}$ such that $w(\mathcal{S}) \leq \alpha \cdot w(\Gamma)$ for some constant $\alpha$. 

By \cref{prop:MS-line-graph}, one can obtain approximation algorithms for contraction sequences from approximation algorithms for tree decompositions; a summary of different approximation algorithms for tree decompositions can be found in \cite{Korhonen:2-approx}. To obtain \cref{co:tractable}, however, we appeal to the polynomial-time algorithm for finding the optimal carving decomposition on planar diagrams in \cite{ST:ratcatcher}, which naturally gives a contraction sequence in view of tree embeddings (see \cite{Bienstock:tree,OGorman:FPT}), and then use the following form of \cref{eq.w-cw}:
\begin{equation}
  w(\Gamma) \leq \frac{3}{2} cw(\Gamma) \leq \frac{3}{2} w(\Gamma).
\end{equation}
This finishes the proof of \cref{co:tractable}.

%\begin{proposition}[Korhonen]
%  \label{prop:Korh-approx}
%  There is an algorithm that, given an $n$-vertex graph $\Gamma$ and an integer $w$, in time $2^{10.8w}n$ either outputs a tree decomposition of $\Gamma$ of width at most $2w+1$, or determines that the tree-width of $\Gamma$ is larger than $w$.
%\end{proposition}

%Since constructing the line graph $\Gamma^*$ is polynomial in the size of $\Gamma$, by constructing $\Gamma^*$, applying \cref{prop:Korh-approx} to $\Gamma^*$, and then \cref{prop:MS-line-graph} to the tree decomposition returned, we obtain 

%\begin{lemma}
%  \label{lem:approx-contraction}
%  There is an algorithm that, given an $n$-vertex graph $\Gamma$ and an integer $w$, in time $2^{10.8w}n$ either outputs a contraction sequence $\mathcal{S}$ of $\Gamma$ with $w(\mathcal{S}) \leq 2w+1$, or determines that the contraction-width of $\Gamma$ is larger than $w$.
%\end{lemma}

%Therefore \cref{co:tractable} follows by combining \cref{lem:approx-contraction} with \cref{th:tractable}.

\section{Benchmarks}
\label{sec:benchmark}

In addition to the benchmarks shown in \cref{fig:time-plot}, we compare the performance of our implementation with others in this section. Although our implementation works for general Reshetikhin--Turaev knot polynomials, we restrict ourselves to the most extensively studied family, the (colored) Jones polynomial, due to the scarcity of well-implemented publicly available packages for the computation of other families. 

\subsection{Colored Jones polynomials on the Conway knot}

For the colored Jones polynomial, the only two contenders are the \texttt{KnotTheory} \cite{KnotTheory} and Hajij and Levitt's \texttt{ColorJones} package \cite{HL:efficient-cj}, both written in \texttt{Mathematica}. Thus we compare with them with a \texttt{Mathematica} implementation of our algorithm, which is obtained from minor edits to the code in \cite{GS:VnData}, where we use an $O(n^3)$ greedy search to produce the contraction sequence; in general, the implementation of our algorithm in \texttt{SnapPy} is at least a magnitude faster than that in \texttt{Mathematica}, partly due to the difference of programming languages. Since we have already drawn \cref{fig:time-width6} showcasing the computation on the knot $11n34$ in \texttt{SnapPy} with a time constraint of roughly $6\cdot 10^3$ seconds, we do the same here to have a grasp of the cross-language performance difference along the way. The benchmark was performed on a MacBook Air with M3 chip and 16GB RAM. This produces \cref{tab:cj-benchmark}.

\begin{table}[!htbp]
  \centering
    \begin{tabular}{ccccc}
    \toprule
    Coloring & \texttt{KnotTheory} & \texttt{ColorJones} & Ours (\texttt{Mathematica}) & Ours (\texttt{Python})\\
    \midrule
    1 & 0.10552 & 0.512003 & 0.977068 & 0.002832\\
    2 & 5.72895 & 0.122924 & 0.968588 & 0.008507\\
    3 & 198.108 & 0.616438 & 1.18240 & 0.033057\\
    4 & OOT & 4.31470 & 1.84138 & 0.062606\\
    5 & & 22.3089 & 4.18859 & 0.208873\\
    6 & & 97.3698 & 9.60003 & 0.319684\\
    7 & & 406.974 & 23.4769 & 0.942618\\
    8 & & 1805.60 & 50.6798 & 1.19045\\
    9 & & OOT & 109.212 & 3.54350\\
    10 & & & 215.333 & 3.90816\\
    11 & & & 419.588 & 10.2546\\
    12 & & & 766.792 & 9.56505 \\
    13 & & & 1347.75 & 26.3211 \\
    14 & & & 2275.95 & 23.0595 \\
    15 & & & 3856.06 & 62.2023 \\
    16 & & & OOT & 50.9203 \\
    17 & & & & 137.517 \\
    18 & & & & 104.946 \\
    19 & & & & 267.455 \\
    20 & & & & 202.984 \\
    21 & & & & 500.720 \\
    22 & & & & 368.573 \\
    23 & & & & 894.917 \\
    24 & & & & 634.879 \\
    25 & & & & 1440.03 \\
    26 & & & & 1023.49 \\
    27 & & & & 2317.58 \\
    28 & & & & 1649.33 \\
    29 & & & & 3704.88 \\
    30 & & & & 2586.93 \\
    31 & & & & 5722.20\\
    32 & & & & 3950.25\\
    33 & & & & OOT \\
    \bottomrule
    \end{tabular}
  \vspace{10pt}
  \caption{Run time (in seconds) of different implementations on computing the colored Jones polynomial on the knot $11n34$ with a time constraint of $6\cdot 10^3$ seconds, where OOT denotes that the time exceeded the time constraint ($6\cdot 10^3$ seconds).}
  \label{tab:cj-benchmark}
\end{table}

Note that in our convention the ordinary Jones polynomial is the $1$-colored Jones polynomial, while in \texttt{ColorJones} it is the $2$-colored Jones polynomial as it corresponds to the $2$-dimensional representation of $U_q(\mathfrak{sl}_2)$, which we shifted in \cref{tab:cj-benchmark} to match our convention. We also remark that \texttt{ColorJones} utilizes $8$ CPU cores during a certain stage of the computation, which cannot be disabled without altering its source code and we did not record the time spent in subkernels by using \texttt{Timing} for the benchmark, while our implementations utilize only $1$ CPU core throughout. On that note, our implementation can be made to utilize parallel computation by adopting the framework of, for example, O'Gorman \cite{OGorman:FPT}.

\subsection{The Jones polynomial on large knot diagrams}
\label{subsec:jones}

For the Jones polynomial, the existing algorithms are of two types, one is to resolve the Kauffman bracket in some order, and the other is fixed-parameter tractable with respect to the tree-width as in \cite{Makowsky:treewidth}. Although various ways of optimization for the former can be made for it to work on small examples (i.e. less than 100 crossings), they are mostly heuristic and scale horrifically. Since we consider knots with at least hundreds of crossings in this section, our only viable opponent is the latter, which is meticulously implemented in \texttt{Regina} \cite{regina}, written in \texttt{C++}, and we compare against it with our \texttt{Python} implementation.

We consider two sources of large diagrams: census knots, tabulated by N. Dunfield and the author \cite{Nathan:exceptional,10census,12census}, with large genus, and a collection of random knots.

For the census knots, we consider the first several ones whose twisted Alexander polynomial, defined by N. Dunfield et al. \cite{Dunfield:twisted}, has a degree greater than 200; their names are as follows:
\begin{equation}
  \label{eq:censusknots}
K9\_24,\ 
K9\_25,\ 
K9\_27,\ 
K9\_32,\ 
K9\_33,\ 
K9\_41,\ 
K9\_43,\ 
K9\_44,\ 
K9\_45.
\end{equation}

Calling \texttt{snappy.Manifold.exterior\_to\_link} \cite{DOR:exterior} directly on the knots in \cref{eq:censusknots} produces diagrams with more than 200 crossings, with which we ran the computation of Jones polynomial on a cluster with a 2.60GHz CPU and a memory constraint of 160GB. 
\texttt{Regina} went out of memory (OOM) for all these diagrams, while our implementation was able to complete all of them within an average time of 206 minutes and memory of 1GB; see \cref{table:censusknots}. See also \cref{appen:jones} for their Jones polynomials.

Unlike for the census knots, widths of generic random knots tend to be much smaller, which allows us to go up in crossing numbers significantly for our benchmark. Our collection of random knots come from two sets of samples, the first contains one per crossing number from 3 to 1000, sampled by N. Dunfield, A. Hirani and M. Obeidin \cite{randomknot}, and the second contains several knots with more than 1000 crossings, generated by \texttt{snappy.random\_link}. Since no two diagrams in our collection have the same crossing number, we use $\mathrm{rd}\_1000$ to denote, for example, the knot in our collection with $1000$ crossings. With this notation, the following are the knots where we perform our benchmark on:
\begin{equation}
  \label{eq:rd}
  \begin{aligned}
  &\mathrm{rd}\_510,\ 
  \mathrm{rd}\_520,\ 
  \mathrm{rd}\_530,\ 
  \mathrm{rd}\_540,\ 
  \mathrm{rd}\_550,\
  \mathrm{rd}\_950,\ 
  \mathrm{rd}\_960,\
  \mathrm{rd}\_970,\ 
  \mathrm{rd}\_980,\ \\
  &\mathrm{rd}\_990,\ 
  \mathrm{rd}\_1000,\
  \mathrm{rd}\_2010,\ 
  \mathrm{rd}\_2012,\ 
  \mathrm{rd}\_2025,\ 
  \mathrm{rd}\_3037,\ 
  \mathrm{rd}\_3070,\ 
  \mathrm{rd}\_4012.
  \end{aligned}
\end{equation}
The diagrams and Jones polynomials for the knots in \cref{eq:rd} are available at \cite{somejones}.

Recall that in our \texttt{SnapPy} implementation, we use \texttt{opt\_einsum} \cite{opteinsum} to find contraction sequences. For large diagrams, \texttt{opt\_einsum} by default uses an $O(n^2)$ greedy algorithm to find the contraction sequence on a given long diagram, and we iterate through every arc to cut open and choose the one that gives the smallest contraction-width. This results in an $O(n^3)$ greedy search,
which dominates the time of the computation when the contraction width is relatively small. Bounding the number of arcs we iterate by a small integer (here we used 30) gives an $O(n^2)$ greedy search, with which finding contraction sequence costs little time but we usually end up with a larger width. For example, for $\mathrm{rd}\_3070$, the $O(n^3)$ greedy search spent 1033.05 minutes to produce a contraction sequence with width 23, while the $O(n^2)$ greedy search spent 9.20 minutes and produced a sequence with width 26; using the contraction sequence with width 23, the calculation of tensor contractions took 283.57 minutes, while it took 1015.09 minutes using the one with width 26. Since this nuanced trade-off impacts the performance, we list both of them in our benchmark in \cref{table:randomknots}; we note that we out-perform \texttt{Regina} on most of these diagrams with either of the two approaches.

As a final remark we note that in practice, the time and RAM cost do not solely depend on the contraction-width and may vary as the contraction sequence changes, see for example $\mathrm{rd}\_4012$. This is because the sparseness of the $R$-matrix is heavily exploited in the implementation, hence the cost may vary as the density of intermediate tensors varies, which the contraction-width does not have control over.

%% Bibliography
%%
%% Please use bibtex, 

\bibliography{references}

@article{opteinsum, doi = {10.21105/joss.00753}, url = {https://doi.org/10.21105/joss.00753}, year = {2018}, publisher = {The Open Journal}, volume = {3}, number = {26}, pages = {753}, author = {Smith, Daniel G. a. and Gray, Johnnie}, title = {opt\_einsum - A Python package for optimizing contraction order for einsum-like expressions}, journal = {Journal of Open Source Software} }

@article {DPSS:treewidth,
    AUTHOR = {de Mesmay, Arnaud and Purcell, Jessica and Schleimer, Saul and
              Sedgwick, Eric},
     TITLE = {On the tree-width of knot diagrams},
   JOURNAL = {J. Comput. Geom.},
  FJOURNAL = {Journal of Computational Geometry},
    VOLUME = {10},
      YEAR = {2019},
    NUMBER = {1},
     PAGES = {164--180},
      ISSN = {1920-180X},
   MRCLASS = {57M25},
  MRNUMBER = {3957223},
MRREVIEWER = {Yuanan\ Diao},
       DOI = {10.20382/jocg.v10i1a6},
       URL = {https://doi.org/10.20382/jocg.v10i1a6},
}

@misc{regina,
    author = {Benjamin A. Burton and Ryan Budney and William Pettersson and others},
    title = {Regina: Software for low-dimensional topology},
    howpublished = {{\tt http://\allowbreak regina-normal.\allowbreak github.\allowbreak io/}},
    year = {1999--2026}}

@misc{GL:FK3,
      title={The $\mathrm{FK}_3$ knot polynomial}, 
      author={Stavros Garoufalidis and Shana Yunsheng Li},
      year={2026},
      eprint={2609.01506},
      archivePrefix={arXiv},
      primaryClass={math.GT},
      url={https://arxiv.org/abs/2609.01506}, 
}

@misc{somejones,
  author = {Shana Yunsheng Li},
  title = {some\_jones},
  year = {2026},
  publisher = {GitHub},
  journal = {GitHub repository},
  howpublished = {\url{https://github.com/Shakugannotorch/some_jones}},
}

@misc{randomknot,
  author = {Dunfield, Nathan and Hirani, Anil and Obeidin, Malik},
  title = {random\_knots},
  year = {2021},
  publisher = {GitHub},
  journal = {GitHub repository},
  howpublished = {\url{https://github.com/NathanDunfield/random_knots}},
}

@article {CDW:simplest,
    AUTHOR = {Callahan, Patrick J. and Dean, John C. and Weeks, Jeffrey R.},
     TITLE = {The simplest hyperbolic knots},
   JOURNAL = {J. Knot Theory Ramifications},
  FJOURNAL = {Journal of Knot Theory and its Ramifications},
    VOLUME = {8},
      YEAR = {1999},
    NUMBER = {3},
     PAGES = {279--297},
      ISSN = {0218-2165,1793-6527},
   MRCLASS = {57M25 (57M50)},
  MRNUMBER = {1691433},
MRREVIEWER = {Colin\ C.\ Adams},
       DOI = {10.1142/S0218216599000195},
       URL = {https://doi.org/10.1142/S0218216599000195},
}

@article {CKP:nextsimplest,
    AUTHOR = {Champanerkar, Abhijit and Kofman, Ilya and Patterson, Eric},
     TITLE = {The next simplest hyperbolic knots},
   JOURNAL = {J. Knot Theory Ramifications},
  FJOURNAL = {Journal of Knot Theory and its Ramifications},
    VOLUME = {13},
      YEAR = {2004},
    NUMBER = {7},
     PAGES = {965--987},
      ISSN = {0218-2165,1793-6527},
   MRCLASS = {57M25},
  MRNUMBER = {2101238},
MRREVIEWER = {Carlo\ Petronio},
       DOI = {10.1142/S021821650400355X},
       URL = {https://doi.org/10.1142/S021821650400355X},
}

@article {CKM:500simplest,
    AUTHOR = {Champanerkar, Abhijit and Kofman, Ilya and Mullen, Timothy},
     TITLE = {The 500 simplest hyperbolic knots},
   JOURNAL = {J. Knot Theory Ramifications},
  FJOURNAL = {Journal of Knot Theory and its Ramifications},
    VOLUME = {23},
      YEAR = {2014},
    NUMBER = {12},
     PAGES = {1450055, 34},
      ISSN = {0218-2165,1793-6527},
   MRCLASS = {57M25 (57M27)},
  MRNUMBER = {3298204},
MRREVIEWER = {Arunima\ Ray},
       DOI = {10.1142/S0218216514500552},
       URL = {https://doi.org/10.1142/S0218216514500552},
}

@article {DOR:exterior,
    AUTHOR = {Dunfield, Nathan M. and Obeidin, Malik and Rudd, Cameron
              Gates},
     TITLE = {Computing a link diagram from its exterior},
   JOURNAL = {Discrete Comput. Geom.},
  FJOURNAL = {Discrete \& Computational Geometry. An International Journal
              of Mathematics and Computer Science},
    VOLUME = {71},
      YEAR = {2024},
    NUMBER = {1},
     PAGES = {121--159},
      ISSN = {0179-5376,1432-0444},
   MRCLASS = {57K10 (57K30)},
  MRNUMBER = {4685711},
       DOI = {10.1007/s00454-023-00533-w},
       URL = {https://doi.org/10.1007/s00454-023-00533-w},
}

@incollection {Nathan:exceptional,
    AUTHOR = {Dunfield, Nathan M.},
     TITLE = {A census of exceptional {D}ehn fillings},
 BOOKTITLE = {Characters in low-dimensional topology},
    SERIES = {Contemp. Math.},
    VOLUME = {760},
     PAGES = {143--155},
 PUBLISHER = {Amer. Math. Soc., [Providence], RI},
      YEAR = {[2020] \copyright 2020},
      ISBN = {978-1-4704-5209-4; 978-1-4704-6135-5},
   MRCLASS = {57K32 (57K35)},
  MRNUMBER = {4193924},
MRREVIEWER = {Scott\ A.\ Taylor},
       DOI = {10.1090/conm/760/15289},
       URL = {https://doi.org/10.1090/conm/760/15289},
}

@article {ST:ratcatcher,
    AUTHOR = {Seymour, P. D. and Thomas, R.},
     TITLE = {Call routing and the ratcatcher},
   JOURNAL = {Combinatorica},
  FJOURNAL = {Combinatorica. An International Journal on Combinatorics and
              the Theory of Computing},
    VOLUME = {14},
      YEAR = {1994},
    NUMBER = {2},
     PAGES = {217--241},
      ISSN = {0209-9683},
   MRCLASS = {05C90 (90B12)},
  MRNUMBER = {1289074},
MRREVIEWER = {Klaus\ Jansen},
       DOI = {10.1007/BF01215352},
       URL = {https://doi.org/10.1007/BF01215352},
}

@InProceedings{OGorman:FPT,
  author =	{O'Gorman, Bryan},
  title =	{{Parameterization of Tensor Network Contraction}},
  booktitle =	{14th Conference on the Theory of Quantum Computation, Communication and Cryptography (TQC 2019)},
  pages =	{10:1--10:19},
  series =	{Leibniz International Proceedings in Informatics (LIPIcs)},
  ISBN =	{978-3-95977-112-2},
  ISSN =	{1868-8969},
  year =	{2019},
  volume =	{135},
  editor =	{van Dam, Wim and Man\v{c}inska, Laura},
  publisher =	{Schloss Dagstuhl -- Leibniz-Zentrum f{\"u}r Informatik},
  address =	{Dagstuhl, Germany},
  URL =		{https://drops.dagstuhl.de/entities/document/10.4230/LIPIcs.TQC.2019.10},
  URN =		{urn:nbn:de:0030-drops-104025},
  doi =		{10.4230/LIPIcs.TQC.2019.10}
}

@misc{HL:efficient-cj,
      title={An Efficient Algorithm to Compute the Colored Jones Polynomial}, 
      author={Mustafa Hajij and Jesse Levitt},
      year={2018},
      eprint={1804.07910},
      archivePrefix={arXiv},
      primaryClass={math.QA},
      url={https://arxiv.org/abs/1804.07910}, 
}

@incollection {KR:colored-jones,
    AUTHOR = {Kirillov, A. N. and Reshetikhin, N. Yu.},
     TITLE = {Representations of the algebra {${U}_q({\rm
              sl}(2)),\;q$}-orthogonal polynomials and invariants of links},
 BOOKTITLE = {Infinite-dimensional {L}ie algebras and groups
              ({L}uminy-{M}arseille, 1988)},
    SERIES = {Adv. Ser. Math. Phys.},
    VOLUME = {7},
     PAGES = {285--339},
 PUBLISHER = {World Sci. Publ., Teaneck, NJ},
      YEAR = {1989},
      ISBN = {9971-50-928-8},
   MRCLASS = {17B35 (05A30 16A24 33A65 57M25 81C40)},
  MRNUMBER = {1026957},
MRREVIEWER = {Bruce\ W.\ Westbury},
}

@article{Orus:tensor,
	title = {Tensor networks for complex quantum systems},
	volume = {1},
	issn = {2522-5820},
	url = {https://doi.org/10.1038/s42254-019-0086-7},
	doi = {10.1038/s42254-019-0086-7},
	number = {9},
	journal = {Nature Reviews Physics},
	author = {Or{\'u}s, Rom{\'a}n},
	month = sep,
	year = {2019},
	pages = {538--550},
}

@article {Bienstock:tree,
    AUTHOR = {Bienstock, Dan},
     TITLE = {On embedding graphs in trees},
   JOURNAL = {J. Combin. Theory Ser. B},
  FJOURNAL = {Journal of Combinatorial Theory. Series B},
    VOLUME = {49},
      YEAR = {1990},
    NUMBER = {1},
     PAGES = {103--136},
      ISSN = {0095-8956,1096-0902},
   MRCLASS = {05C10 (05C05)},
  MRNUMBER = {1056822},
MRREVIEWER = {Miros\l aw\ Truszczy\'nski},
       DOI = {10.1016/0095-8956(90)90066-9},
       URL = {https://doi.org/10.1016/0095-8956(90)90066-9},
}

@article {KS:TPC,
    AUTHOR = {Korach, Ephraim and Solel, Nir},
     TITLE = {Tree-width, path-width, and cutwidth},
   JOURNAL = {Discrete Appl. Math.},
  FJOURNAL = {Discrete Applied Mathematics. The Journal of Combinatorial
              Algorithms, Informatics and Computational Sciences},
    VOLUME = {43},
      YEAR = {1993},
    NUMBER = {1},
     PAGES = {97--101},
      ISSN = {0166-218X,1872-6771},
   MRCLASS = {05C85 (68Q25 68R10)},
  MRNUMBER = {1218045},
       DOI = {10.1016/0166-218X(93)90171-J},
       URL = {https://doi.org/10.1016/0166-218X(93)90171-J},
}

@book {ModernComputerAlgebra,
    AUTHOR = {von zur Gathen, Joachim and Gerhard, J\"urgen},
     TITLE = {Modern computer algebra},
 PUBLISHER = {Cambridge University Press, New York},
      YEAR = {1999},
     PAGES = {xiv+753},
      ISBN = {0-521-64176-4},
   MRCLASS = {68W30 (11Y16 68-01 68-02)},
  MRNUMBER = {1689167},
MRREVIEWER = {Jeffrey\ O.\ Shallit},
}

@incollection {Korhonen:2-approx,
    AUTHOR = {Korhonen, Tuukka},
     TITLE = {A single-exponential time 2-approximation algorithm for
              treewidth},
 BOOKTITLE = {2021 {IEEE} 62nd {A}nnual {S}ymposium on {F}oundations of
              {C}omputer {S}cience---{FOCS} 2021},
     PAGES = {184--192},
 PUBLISHER = {IEEE Computer Soc., Los Alamitos, CA},
      YEAR = {[2022] \copyright 2022},
      ISBN = {978-1-6654-2055-6},
   MRCLASS = {68W25},
  MRNUMBER = {4399680},
       DOI = {10.1109/FOCS52979.2021.00026},
       URL = {https://doi.org/10.1109/FOCS52979.2021.00026},
}

@incollection {Maria:FPT,
    AUTHOR = {Maria, Cl\'ement},
     TITLE = {Parameterized complexity of quantum knot invariants},
 BOOKTITLE = {37th {I}nternational {S}ymposium on {C}omputational
              {G}eometry},
    SERIES = {LIPIcs. Leibniz Int. Proc. Inform.},
    VOLUME = {189},
     PAGES = {Art. No. 53, 17},
 PUBLISHER = {Schloss Dagstuhl. Leibniz-Zent. Inform., Wadern},
      YEAR = {2021},
      ISBN = {978-3-95977-184-9},
   MRCLASS = {68Q27 (57K16)},
  MRNUMBER = {4287060},
}

@article {Bar-Natan:Khovanov,
    AUTHOR = {Bar-Natan, Dror},
     TITLE = {Khovanov's homology for tangles and cobordisms},
   JOURNAL = {Geom. Topol.},
  FJOURNAL = {Geometry and Topology},
    VOLUME = {9},
      YEAR = {2005},
     PAGES = {1443--1499},
      ISSN = {1465-3060,1364-0380},
   MRCLASS = {57M27 (57M25 57R56)},
  MRNUMBER = {2174270},
MRREVIEWER = {Justin\ Sawon},
       DOI = {10.2140/gt.2005.9.1443},
       URL = {https://doi.org/10.2140/gt.2005.9.1443},
}

@book {Ohtsuki:quantum,
    AUTHOR = {Ohtsuki, Tomotada},
     TITLE = {Quantum invariants},
    SERIES = {Series on Knots and Everything},
    VOLUME = {29},
      NOTE = {A study of knots, 3-manifolds, and their sets},
 PUBLISHER = {World Scientific Publishing Co., Inc., River Edge, NJ},
      YEAR = {2002},
     PAGES = {xiv+489},
      ISBN = {981-02-4675-7},
   MRCLASS = {57M27},
  MRNUMBER = {1881401},
MRREVIEWER = {Justin\ D.\ Roberts},
}

@article {LT:separator,
    AUTHOR = {Lipton, Richard J. and Tarjan, Robert Endre},
     TITLE = {A separator theorem for planar graphs},
   JOURNAL = {SIAM J. Appl. Math.},
  FJOURNAL = {SIAM Journal on Applied Mathematics},
    VOLUME = {36},
      YEAR = {1979},
    NUMBER = {2},
     PAGES = {177--189},
      ISSN = {0036-1399},
   MRCLASS = {68E10},
  MRNUMBER = {524495},
MRREVIEWER = {Ian\ Munro},
       DOI = {10.1137/0136016},
       URL = {https://doi.org/10.1137/0136016},
}

@book {CFKLMPPS:parametrized,
    AUTHOR = {Cygan, Marek and Fomin, Fedor V. and Kowalik, \L ukasz and
              Lokshtanov, Daniel and Marx, D\'aniel and Pilipczuk, Marcin
              and Pilipczuk, Micha\l{} and Saurabh, Saket},
     TITLE = {Parameterized algorithms},
 PUBLISHER = {Springer, Cham},
      YEAR = {2015},
     PAGES = {xviii+613},
      ISBN = {978-3-319-21274-6; 978-3-319-21275-3},
   MRCLASS = {68-01 (05-01 05C85 68Q17 68Q25 68R10 68Wxx)},
  MRNUMBER = {3380745},
MRREVIEWER = {Henning\ Fernau},
       DOI = {10.1007/978-3-319-21275-3},
       URL = {https://doi.org/10.1007/978-3-319-21275-3},
}

@article {MS:tensor,
    AUTHOR = {Markov, Igor L. and Shi, Yaoyun},
     TITLE = {Simulating quantum computation by contracting tensor networks},
   JOURNAL = {SIAM J. Comput.},
  FJOURNAL = {SIAM Journal on Computing},
    VOLUME = {38},
      YEAR = {2008},
    NUMBER = {3},
     PAGES = {963--981},
      ISSN = {0097-5397,1095-7111},
   MRCLASS = {81P68 (05C83 68Q05 68Q10 68R10)},
  MRNUMBER = {2421074},
MRREVIEWER = {Hari\ Kiran\ Krovi},
       DOI = {10.1137/050644756},
       URL = {https://doi.org/10.1137/050644756},
}

@article {MK:Jones,
    AUTHOR = {Meichanetzidis, Konstantinos and Kourtis, Stefanos},
     TITLE = {Evaluating the {J}ones polynomial with tensor networks},
   JOURNAL = {Phys. Rev. E},
  FJOURNAL = {Physical Review E},
    VOLUME = {100},
      YEAR = {2019},
    NUMBER = {3},
     PAGES = {033303, 7},
      ISSN = {2470-0045,2470-0053},
   MRCLASS = {57M27},
  MRNUMBER = {4019495},
       DOI = {10.1103/physreve.100.033303},
       URL = {https://doi.org/10.1103/physreve.100.033303},
}

@article{Makowsky:treewidth,
title = {Coloured Tutte polynomials and Kauffman brackets for graphs of bounded tree width},
journal = {Discrete Applied Mathematics},
volume = {145},
number = {2},
pages = {276-290},
year = {2005},
note = {Structural Decompositions, Width Parameters, and Graph Labelings},
issn = {0166-218X},
doi = {https://doi.org/10.1016/j.dam.2004.01.016},
url = {https://www.sciencedirect.com/science/article/pii/S0166218X04002513},
author = {J.A. Makowsky}
}

@article {Turaev:YB,
    AUTHOR = {Turaev, V. G.},
     TITLE = {The {Y}ang-{B}axter equation and invariants of links},
   JOURNAL = {Invent. Math.},
  FJOURNAL = {Inventiones Mathematicae},
    VOLUME = {92},
      YEAR = {1988},
    NUMBER = {3},
     PAGES = {527--553},
      ISSN = {0020-9910,1432-1297},
   MRCLASS = {57M25 (20F36 82A68)},
  MRNUMBER = {939474},
MRREVIEWER = {Toshitake\ Kohno},
       DOI = {10.1007/BF01393746},
       URL = {https://doi.org/10.1007/BF01393746},
}

@article {Shumakovitch:odd,
    AUTHOR = {Shumakovitch, Alexander N.},
     TITLE = {Patterns in odd {K}hovanov homology},
   JOURNAL = {J. Knot Theory Ramifications},
  FJOURNAL = {Journal of Knot Theory and its Ramifications},
    VOLUME = {20},
      YEAR = {2011},
    NUMBER = {1},
     PAGES = {203--222},
      ISSN = {0218-2165,1793-6527},
   MRCLASS = {57M27 (57R17)},
  MRNUMBER = {2777025},
MRREVIEWER = {Radmila\ Sazdanovi\'c},
       DOI = {10.1142/S0218216511008802},
       URL = {https://doi.org/10.1142/S0218216511008802},
}

@article {JVW:complexity,
    AUTHOR = {Jaeger, F. and Vertigan, D. L. and Welsh, D. J. A.},
     TITLE = {On the computational complexity of the {J}ones and {T}utte
              polynomials},
   JOURNAL = {Math. Proc. Cambridge Philos. Soc.},
  FJOURNAL = {Mathematical Proceedings of the Cambridge Philosophical
              Society},
    VOLUME = {108},
      YEAR = {1990},
    NUMBER = {1},
     PAGES = {35--53},
      ISSN = {0305-0041,1469-8064},
   MRCLASS = {05B35 (57M25 68Q25)},
  MRNUMBER = {1049758},
MRREVIEWER = {Mark\ E.\ Kidwell},
       DOI = {10.1017/S0305004100068936},
       URL = {https://doi.org/10.1017/S0305004100068936},
}

@incollection {Witten,
    AUTHOR = {Witten, Edward},
     TITLE = {Quantum field theory and the {J}ones polynomial},
 BOOKTITLE = {Braid group, knot theory and statistical mechanics},
    SERIES = {Adv. Ser. Math. Phys.},
    VOLUME = {9},
     PAGES = {239--329},
 PUBLISHER = {World Sci. Publ., Teaneck, NJ},
      YEAR = {1989},
      ISBN = {9971-50-828-1; 9971-50-833-8},
   MRCLASS = {57M25 (58F07 58G40 81R50 81T30 81T40)},
  MRNUMBER = {1062429},
       DOI = {10.1142/9789812798350\_0009},
       URL = {https://doi.org/10.1142/9789812798350_0009},
}

@article {GL:patterns,
    title={Patterns of the ${V}_2$-polynomial of knots}, 
      author={Stavros Garoufalidis and Shana Yunsheng Li},
      journal = {Experimental Mathematics},
      year={2026},
      pages = {1-31},
      doi = {10.1080/10586458.2026.2651081}
}

@misc{GLY:LG,
  author = {Stavros Garoufalidis and Shana Yunsheng Li and Josephine Yu},
  title = {Positivity and concavity of the colored {L}inks--{G}ould polynomials of knots},
  eprint={2605.30754},
  archivePrefix={arXiv},
  primaryClass={math.GT},
  url={https://arxiv.org/abs/2605.30754}, 
}

@article {BNV:API,
    AUTHOR = {Bar-Natan, Dror and van der Veen, Roland},
     TITLE = {A perturbed-{A}lexander invariant},
   JOURNAL = {Quantum Topol.},
  FJOURNAL = {Quantum Topology},
    VOLUME = {15},
      YEAR = {2024},
    NUMBER = {3},
     PAGES = {449--472},
      ISSN = {1663-487X},
   MRCLASS = {57K14 (16T99)},
  MRNUMBER = {4814687},
MRREVIEWER = {Zhiqing Yang},
       DOI = {10.4171/qt/206},
       URL = {https://doi.org/10.4171/qt/206},
}

@article {JC:tensornet,
    AUTHOR = {Bridgeman, Jacob and Chubb, Christopher},
     TITLE = {Hand-waving and interpretive dance: an introductory course on
              tensor networks},
   JOURNAL = {J. Phys. A},
  FJOURNAL = {Journal of Physics. A. Mathematical and Theoretical},
    VOLUME = {50},
      YEAR = {2017},
    NUMBER = {22},
     PAGES = {223001, 61},
      ISSN = {1751-8113},
   MRCLASS = {81P45},
  MRNUMBER = {3659106},
       DOI = {10.1088/1751-8121/aa6dc3},
       URL = {https://doi.org/10.1088/1751-8121/aa6dc3},
}

@misc{snappy,
     author={Culler, Marc and Dunfield, Nathan M. and Goerner,
     Matthias and Weeks, Jeffrey R.},
     title={Snap{P}y, a computer program for studying the geometry and topology of $3$-manifolds},
     howpublished={Available at \url{http://snappy.computop.org} (26/10/2026)},
}

@unpublished{12census,
  author = {Shana Yunsheng Li},
  title  = {The complete 11- and 12-tetrahedra census of orientable cusped hyperbolic 3-manifolds},
  year   = {2026},
  note   = {In preparation}
}

@incollection {10census,
    AUTHOR = {Li, Shana Yunsheng},
     TITLE = {The complete 10-tetrahedra census of orientable cusped
              hyperbolic 3-manifolds},
 BOOKTITLE = {42nd {I}nternational {S}ymposium on {C}omputational
              {G}eometry},
    SERIES = {LIPIcs. Leibniz Int. Proc. Inform.},
    VOLUME = {367},
     PAGES = {Paper No. 73, 15},
 PUBLISHER = {Schloss Dagstuhl. Leibniz-Zent. Inform., Wadern},
      YEAR = {2026},
      ISBN = {978-3-95977-418-5},
   MRCLASS = {57K32 (57-08 57-11)},
  MRNUMBER = {5080744},
       DOI = {10.4230/lipics.socg.2026.73},
       URL = {https://doi.org/10.4230/lipics.socg.2026.73},
}

@Misc {KnotTheory,
    AUTHOR = {},
     TITLE = {The {K}not {A}tlas},
HOWPUBLISHED = {\url{http://katlas.org/}},
urldate  = {2026-10-19},
}

@Misc {spherogram,
    AUTHOR = {Culler, Marc and Dunfield, Nathan and Lipshitz, Robert and Dickinson, Jennet and Obeidin, Malik},
     TITLE = {Spherogram},
     year = {2026},
     publisher = {GitHub},
     journal = {GitHub repository},
HOWPUBLISHED = {\url{https://github.com/3-manifolds/Spherogram}},
}

@article {Dunfield:twisted,
    AUTHOR = {Dunfield, Nathan and Friedl, Stefan and Jackson, Nicholas},
     TITLE = {Twisted {A}lexander polynomials of hyperbolic knots},
   JOURNAL = {Exp. Math.},
  FJOURNAL = {Experimental Mathematics},
    VOLUME = {21},
      YEAR = {2012},
    NUMBER = {4},
     PAGES = {329--352},
      ISSN = {1058-6458},
   MRCLASS = {57M25},
  MRNUMBER = {3004250},
       DOI = {10.1080/10586458.2012.669268},
       URL = {https://doi.org/10.1080/10586458.2012.669268},
}

@article {GK:multi,
    AUTHOR = {Garoufalidis, Stavros and Kashaev, Rinat},
     TITLE = {Multivariable {K}not {P}olynomials from {B}raided {H}opf
              {A}lgebras with {A}utomorphisms},
   JOURNAL = {Publ. Res. Inst. Math. Sci.},
  FJOURNAL = {Publications of the Research Institute for Mathematical
              Sciences},
    VOLUME = {62},
      YEAR = {2026},
    NUMBER = {1},
     PAGES = {75--114},
      ISSN = {0034-5318,1663-4926},
   MRCLASS = {57K14 (16 57K10)},
  MRNUMBER = {5032319},
       DOI = {10.4171/prims/62-1-3},
       URL = {https://doi.org/10.4171/prims/62-1-3},
}

@Misc {Vn,
    AUTHOR = {Garoufalidis, Stavros and Matthew Harper and Kashaev, Rinat and
             Kohli, Ben-Michael and Wagner, Emmanuel},
     TITLE = {On the colored {L}inks--{G}ould polynomial},
	NOTE = {Preprint 2025,
             \href{https://arxiv.org/abs/2509.10911}{arXiv:2509.10911}},
}

@data {GS:VnData,
    AUTHOR = {Garoufalidis, Stavros and Li, Shana},
     TITLE = {Values of {$V_n$}-polynomials of knots},
 PUBLISHER = {Harvard Dataverse},
       UNF = {UNF:6:xsUpf4nwo557+RqL7sHv9w==},
      YEAR = {2024},
   VERSION = {V1},
       DOI = {10.7910/DVN/XE4TOF},
}

@article {HTW,
    AUTHOR = {Hoste, Jim and Thistlethwaite, Morwen and Weeks, Jeff},
     TITLE = {The first 1,701,936 knots},
   JOURNAL = {Math. Intelligencer},
  FJOURNAL = {The Mathematical Intelligencer},
    VOLUME = {20},
      YEAR = {1998},
    NUMBER = {4},
     PAGES = {33--48},
      ISSN = {0343-6993},
   MRCLASS = {57M25},
  MRNUMBER = {1646740},
MRREVIEWER = {C. Kearton},
       DOI = {10.1007/BF03025227},
       URL = {https://doi.org/10.1007/BF03025227},
}

@article {Jones,
    AUTHOR = {Jones, Vaughan},
     TITLE = {Hecke algebra representations of braid groups and link
              polynomials},
   JOURNAL = {Ann. of Math. (2)},
  FJOURNAL = {Annals of Mathematics. Second Series},
    VOLUME = {126},
      YEAR = {1987},
    NUMBER = {2},
     PAGES = {335--388},
      ISSN = {0003-486X},
     CODEN = {ANMAAH},
   MRCLASS = {46L99 (20F36 22D25 46L35 46L55 57M25)},
  MRNUMBER = {MR908150 (89c:46092)},
MRREVIEWER = {Pierre de la Harpe},
       DOI = {10.2307/1971403},
}

@article{RJF:optimaltensor,
author = {Pfeifer, Robert and Haegeman, Jutho and Verstraete, Frank},
address = {United States},
issn = {1539-3755},
journal = {Physical Review E},
language = {eng},
number = {3},
pages = {033315-},
publisher = {American Physical Society (APS)},
title = {Faster identification of optimal contraction sequences for tensor networks},
volume = {90},
year = {2014},
}

@article {RT:ribbon,
    AUTHOR = {Reshetikhin, Nikolai and Turaev, Vladimir},
     TITLE = {Ribbon graphs and their invariants derived from quantum
              groups},
   JOURNAL = {Comm. Math. Phys.},
  FJOURNAL = {Communications in Mathematical Physics},
    VOLUME = {127},
      YEAR = {1990},
    NUMBER = {1},
     PAGES = {1--26},
      ISSN = {0010-3616},
   MRCLASS = {57M25 (16W30 17B35)},
  MRNUMBER = {1036112},
MRREVIEWER = {Louis H. Kauffman},
       URL = {http://projecteuclid.org/euclid.cmp/1104180037},
}

\section{AI Statement}

The implementation of tensor contractions and Laurent polynomials in \texttt{SnapPy} was produced and optimized by Claude Code Opus under thorough human supervision from the author.

\appendix

\section{The Jones polynomial of census knots in \cref{table:censusknots}}
\label{appen:jones}

We present the Jones polynomials of the census knots in \cref{table:censusknots} in this appendix, obtained from direct computation by our implementation in \texttt{SnapPy}. We remark that their Jones polynomials are of especially simple forms compared to generic knots with similar numbers of crossings, potentially due to them being special kinds of twisted-torus knots as in the spirit of \cite{CDW:simplest,CKP:nextsimplest,CKM:500simplest}. In comparison, the Jones polynomial of $\textrm{rd}\_210$ in \cite{randomknot} consists of 152 monomials with the largest coefficient being $\numprint{2809888906534799157919148179}$. 

\begin{equation*}
  \begin{aligned}
  J_{K9\_24}(q) =
  & q^{-177} - 2 q^{-176} + 2 q^{-175} - 2 q^{-174} + 2 q^{-173} - 2 q^{-172} + 2 q^{-171} - 2 q^{-170}\\
  & + 2 q^{-169} - 2 q^{-168} + 2 q^{-167} - 2 q^{-166} + q^{-165} - q^{-164} + q^{-163} - q^{-162} \\
  &+ q^{-161} + q^{-152} - q^{-151} + q^{-150} - q^{-149} + q^{-148} - q^{-147} + q^{-146} \\
  &- 2 q^{-145}+ 2 q^{-144} - 2 q^{-143} + 2 q^{-142} - 2 q^{-141} + 2 q^{-140} - 2 q^{-139} \\
  &+ 2 q^{-138} - 2 q^{-137} + 2 q^{-136} - 2 q^{-135} + q^{-134} - 2 q^{-133} + q^{-132} - 2 q^{-131} \\
  & + q^{-130} - 2 q^{-129} + q^{-128} - 2 q^{-127} + q^{-126} - q^{-125} + q^{-124} - q^{-123} \\
  &+ q^{-122} - q^{-121} + q^{-120} - q^{-119} + q^{-118} - q^{-117} + q^{-116} - q^{-115} \\
  &+ q^{-114} - q^{-113} + q^{-112} + q^{-110} + q^{-108} + q^{-106} + q^{-104} + q^{-102}  
  \end{aligned}
\end{equation*}

\begin{equation*}
  \begin{aligned}
  J_{K9\_25}(q) =
  & q^{-200} - q^{-199} + q^{-198} - 2 q^{-197} + 2 q^{-196} - 2 q^{-195} + 2 q^{-194} - 2 q^{-193} \\
  &  + 2 q^{-192} - 2 q^{-191} + q^{-190} - q^{-189} + q^{-188} - q^{-187} + q^{-186} - q^{-153} \\
  &  - q^{-151} - q^{-149} - q^{-147} - q^{-145} + q^{-144} - q^{-143} + q^{-142} - q^{-141} \\
  &  + q^{-140} - q^{-139} + q^{-138} - q^{-137} + q^{-136} - q^{-135} + q^{-134} - q^{-133} \\
  &  + q^{-132} - q^{-131} + q^{-130} - q^{-129} + q^{-128} + q^{-126} + q^{-124} + q^{-122} \\
  &  + q^{-120} + q^{-118}
  \end{aligned}
\end{equation*}

\begin{equation*}
  \begin{aligned}
  J_{K9\_27}(q) =
  & q^{114} + q^{116} + q^{118} + q^{120} + q^{122} + q^{124} - q^{125} + q^{126} \\
&  - q^{127} + q^{128} - q^{129} + q^{130} - q^{131} + q^{132} - q^{133} + q^{134} \\
&  - q^{135} + q^{136} - q^{137} + q^{138} - q^{139} + q^{140} - 2 q^{141} + q^{142} \\
&  - 2 q^{143} + q^{144} - 2 q^{145} + q^{146} - 2 q^{147} + q^{148} - 2 q^{149} + 2 q^{150} \\
&  - 2 q^{151} + 2 q^{152} - 2 q^{153} + 2 q^{154} - 2 q^{155} + 2 q^{156} - 2 q^{157} + 2 q^{158} \\
&  - 2 q^{159} + 2 q^{160} - 2 q^{161} + 2 q^{162} - q^{163} + q^{164} - q^{165} + q^{166} \\
&  - q^{167} + q^{168} - q^{169} - q^{180} + q^{181} - q^{182} + q^{183} - q^{184} \\
&  + 2 q^{185} - 2 q^{186} + 2 q^{187} - 2 q^{188} + 2 q^{189} - 2 q^{190} + 2 q^{191} - 2 q^{192} \\
&  + 2 q^{193} - 2 q^{194} + 2 q^{195} - 2 q^{196} + 2 q^{197} - q^{198}
  \end{aligned}
\end{equation*}

\begin{equation*}
  \begin{aligned}
  J_{K9\_32}(q) =
  & q^{-204} - q^{-203} - q^{-195} + q^{-194} - q^{-193} + q^{-192} + q^{-179} - q^{-178} \\
&  + q^{-177} - q^{-176} + q^{-175} - q^{-174} - q^{-156} - q^{-154} - q^{-152} - q^{-150} \\
&  + q^{-126} + q^{-124} + q^{-122} + q^{-120} + q^{-118}
  \end{aligned}
\end{equation*}

\begin{equation*}
  \begin{aligned}
  J_{K9\_33}(q) =
  & -q^{-217} + q^{-216} + q^{-207} - q^{-206} + q^{-205} - q^{-204} - q^{-190} + q^{-189} \\
&  - q^{-188} + q^{-187} - q^{-186} + q^{-185} - q^{-165} - q^{-163} - q^{-161} - q^{-159} \\
&  + q^{-134} + q^{-132} + q^{-130} + q^{-128} + q^{-126}
  \end{aligned}
\end{equation*}

\begin{equation*}
  \begin{aligned}
  J_{K9\_41}(q) =
 & q^{128} + q^{130} + q^{132} + q^{134} + q^{136} + q^{138} - q^{156} - q^{158} \\
&  - q^{160} - q^{162} - q^{164} - q^{178} + q^{179} - q^{180} + q^{181} - q^{182} \\
&  + q^{183} - q^{184} + q^{185} + q^{196} - q^{197} + q^{198} - q^{199} + q^{200} \\
&  - q^{201} - q^{209} + q^{210} - q^{211} + q^{212} + q^{217} - q^{218}
  \end{aligned}
\end{equation*}

\begin{equation*}
  \begin{aligned}
  J_{K9\_43}(q) =
 & q^{137} + q^{139} + q^{141} + q^{143} + q^{145} + q^{147} - q^{166} - q^{168} \\
&  - q^{170} - q^{172} - q^{174} + q^{190} - q^{191} + q^{192} - q^{193} + q^{194} \\
&  - q^{195} + q^{196} - q^{197} - q^{209} + q^{210} - q^{211} + q^{212} - q^{213} \\
&  + q^{214} + q^{223} - q^{224} + q^{225} - q^{226} - q^{232} + q^{233}
  \end{aligned}
\end{equation*}

\begin{equation*}
  \begin{aligned}
  J_{K9\_44}(q) =
 & q^{126} + q^{128} + q^{130} + q^{132} + q^{134} + q^{136} - q^{153} - q^{155} \\
&  - q^{157} - q^{159} - q^{161} + q^{193} - q^{194} + q^{195} - q^{196} + q^{197} \\
&  - q^{198} - q^{206} + q^{207} - q^{208} + q^{209}
  \end{aligned}
\end{equation*}

\begin{equation*}
  \begin{aligned}
  J_{K9\_45}(q) =
 & -q^{-229} + q^{-228} - q^{-227} + q^{-226} + q^{-217} - q^{-216} + q^{-215} - q^{-214} \\
&  + q^{-213} - q^{-212} - q^{-177} - q^{-175} - q^{-173} - q^{-171} - q^{-169} + q^{-149} \\
&  + q^{-147} + q^{-145} + q^{-143} + q^{-141} + q^{-139}
  \end{aligned}
\end{equation*}

\end{document}